\documentclass[11pt]{article} 
\usepackage{a4wide}
\usepackage{amsmath,amsfonts}

\newtheorem{theorem}{Theorem}

\newtheorem{proposition}{Proposition}
\newtheorem{lemma}{Lemma} 
\newtheorem{claim}{Claim}
\numberwithin{theorem}{section}
\numberwithin{claim}{section}
\numberwithin{equation}{section}
\numberwithin{lemma}{section}
\numberwithin{proposition}{section}
 
\begin{document}

\title{Stability of two soliton collision \\ for nonintegrable gKdV equations
\footnote{
This research was supported in part by the Agence Nationale de la Recherche (ANR ONDENONLIN).
}}
\author{Yvan Martel$^{(1)}$ 
  and  Frank Merle$^{(2)}$ }
\date{(1) Universit\'e de Versailles Saint-Quentin-en-Yvelines,
 Math\'ematiques \\
 45, av. des Etats-Unis,
 78035 Versailles cedex, France\\  
 martel@math.uvsq.fr\\
\quad \\ 
(2)
Universit\'e de Cergy-Pontoise, IHES and CNRS,
Math\'ematiques   \\
2, av. Adolphe Chauvin,
95302 Cergy-Pontoise cedex, France \\
 Frank.Merle@math.u-cergy.fr\\
}
\maketitle
 
\begin{abstract}

We continue our study of the collision of two solitons for the 
subcritical generalized KdV equations
\begin{equation}\label{kdvabs}
    \partial_t u + \partial_x (\partial_x^2 u + f(u))=0.
\end{equation}
Solitons are solutions  o the type
$u(t,x)=Q_{c_0}(x-x_0-c_0 t)$ where $c_0>0$.
In \cite{MMcol1}, mainly devoted to the case 
$f(u)=u^4$,
 we have introduced a new framework to understand the collision of 
 two solitons $Q_{c_1}$, $Q_{c_2}$ for \eqref{kdvabs} in the case
$c_2 \ll c_1$ (or equivalently, $\|Q_{c_2}\|_{H^1}\ll \|Q_{c_1}\|_{H^1}$).
In this paper, we consider the case of a general nonlinearity $f(u)$
for which $Q_{c_1}$, $Q_{c_2}$  are nonlinearly stable. In particular,
 since
$f$ is general and $c_1$ can be large, the results are not   pertubations of the ones
for the power case  in \cite{MMcol1}.

First, we prove that the two solitons survive the collision up to a shift in their
trajectory and up to a small perturbation term whose size is explicitely controlled from above: after the collision,
$u(t)\sim Q_{c_1^+}+Q_{c_2^+} $ where $c_j^+$ is close to $c_j$ $(j=1,2)$.
Then, we  exhibit new exceptional solutions similar to multi-soliton solutions:
  for all $c_1,$ $c_2>0$, $c_2\ll c_1$, there exists a
solution $\varphi(t)$ such that
\begin{equation*}
		\varphi(t,x)=Q_{c_1}(x{-}\rho_1(t)) 
		+ Q_{c_2}(x{-}\rho_2(t)) + \eta(t,x), 
		\text{ for $t\ll -1$,}
\end{equation*}
\begin{equation*}
		\varphi(t,x)=Q_{c_1}(x{-}\rho_1(t))
		+ Q_{c_2}(x{-}\rho_2(t)) + \eta(t,x), 
		\text{ for $t\gg 1$,}
\end{equation*}
where $\rho_j(t)\to c_j$ ($j=1,2$)
and $\eta(t)$   converges to $0$ in a neighborhood of the solitons as $t\to \pm\infty$.

The analysis is splitted in two distinct parts.
For the interaction region, we extend the algebraic tools developed in \cite{MMcol1} for the power case,
by expanding $f(u)$ as a sum of powers plus a perturbation term.
To study the solutions in large time, we rely on previous tools 
on asymptotic stability in \cite{MM1}, \cite{MMT} and \cite{MMnonlinearity}, refined in \cite{MMas1},
\cite{MMas2}.
\end{abstract}

\eject

\section{Introduction}

We consider the generalized Korteweg-de Vries (gKdV) equations:
\begin{equation}\label{kdvf}
    \partial_t u + \partial_x (\partial_x^2 u + f(u))=0, 
    \quad (t,x) \in \mathbb{R}^+ \times \mathbb{R}, \quad u(0)=u_0\in H^1(\mathbb{R}),
\end{equation}
for general $C^{s}$ nonlinearity $f$ 
for which small solitons are stable.
We assume that for $p=2,$ $3$ or $4$,
\begin{equation}\label{surf} 
f(u)=u^p + f_1(u) \quad \text{where $f_1$ is $C^{p+4}$  and }\quad 
\lim_{u\to 0} \left|\frac {f_1(u)}{u^p}\right|=0.
\end{equation}
Remark that if the nonlinearity is of the form $f(u)=a u^p + f_1(u)$, $a>0$,
then we may assume $a=1$ by  considering $\widetilde u(t,x)= a^{\frac 1{p-1}}u(t,x)$
instead of $u(t,x)$ and changing $f_1$ accordingly.
We only consider the case where $p=2$, $3$ or $4$ in \eqref{surf}
since otherwise solitons with small speed would  not be stable, which is necessary in this paper.
Denote $F(s)=\int_0^s f(s') ds'$.

The Cauchy problem for equation \eqref{kdvf} is locally 
well-posed in $H^1(\mathbb{R})$ (see Kenig, Ponce and Vega \cite{KPV}).
All solutions considered in this paper are global in time.
For $H^1$ solutions, the following quantities are conserved:
\begin{equation}\label{masse}\begin{split}
&  \int  u^2(t,x)dx=\int  u^2_0(x)dx, \\
&   E(u(t))={\frac 1 2} \int  (\partial_x u)^2(t,x)dx - \int  F(u(t,x))dx
         ={\frac1 2} \int  (\partial_x u_0)^2(x)dx  -  \int F(u_0(x))dx.
\end{split}\end{equation}
Recall that equation \eqref{kdvf} has soliton solutions, i.e. of the
form $u(t,x)=Q_c(x-x_0-c\,t)$ where $c>0$, $x_0\in \mathbb{R}$ and 
\begin{equation}\label{ellipticf}
Q_c''+f(Q_c)=c\,Q_c,\quad Q_c\in H^1.
\end{equation}
 Note that, for all $c>0$,
 if $p=2$, $4$ then there is at most one
solution of \eqref{ellipticf} (up to translations), which is positive,
whereas for $p=3$, it might   exist a positive and a
negative solution of \eqref{ellipticf}.
For all $c>0$, if   a solution $Q_c>0$ of \eqref{ellipticf} exists then it can be chosen even on $\mathbb{R}$ and decreasing on $\mathbb{R}^+$
(and similarly if $Q_c<0$).
We refer to section 6 of Berestycki and Lions \cite{BL} for these properties and a necessary and sufficient condition for existence.

In this paper, we consider only nonlinearly stable solitons in the sense of Weinstein
\cite{We2}, i.e. such that
\begin{equation}\label{stable}
\frac d{dc'} \int {Q_{c'}^2(x) dx }_{\big |{c'}=c}>0.
\end{equation}
Note that since $p=2$, $3$ or $4$ in \eqref{surf}, this condition is
satisfied for $c>0$ small enough. We recall the following stability result.

\medskip

\noindent\textbf{Stability result \cite{We2}.}\textit{
Let $c>0$ be such that \eqref{stable} holds. Then, there exist $K, \alpha_0>0$ such that
for any $u_0\in H^1$, if $\|u_0-Q_c\|_{H^1}\leq \alpha_0$, then 
the solution $u(t)$ of \eqref{kdvf} is global and, for all $t\in \mathbb{R}$,
$\inf_{y\in \mathbb{R}} \|u(t,.+y)-Q_c\|_{H^1}\leq K \alpha_0$.}

\medskip

From \cite{BL} and  
\eqref{surf}, it follows that there exists $c_*(f)>0$ (possibly $+\infty$) defined by
$$
c_*(f)=\sup\{c>0\text{ such that $\forall c'\in (0,c), \exists Q_{c'}$ positive
solution of \eqref{ellipticf}} \}.
$$
In \cite{MMas1}, we have proved that $0<c<c_*(f)$ is a sufficient condition of
asymptotic stability in the energy space $H^1$ around the soliton $Q_c$. Combining the stability result and
the asymptotic stability result, we obtain the following.

\medskip

\noindent\textbf{Asymptotic stability \cite{MM1}, \cite{MMas1}.}\textit{
Let $0<c<c_*(f)$ be such that \eqref{stable} holds.
There exists $\alpha_0>0$ such that  
for any $u_0\in H^1$, if $\|u_0-Q_c\|_{H^1}\leq \alpha_0$, then the solution
$u(t)$ of \eqref{kdvf} is global and  there exist  $c^+\in (0,c_*(f))$, $t\mapsto \rho(t)\in \mathbb{R}$ such that
for all $A>0$,}
\begin{equation}
\label{th1-1}
u(t)-Q_{c^+}(.-\rho(t))\to 0 \quad \text{in $H^1(x>\tfrac c{10} t)$ as $t\to +\infty$.}
\end{equation}

\medskip

We also recall from \cite{Martel} the following result of existence and uniqueness of asymptotic $N$-soliton-like
solutions (see Theorem 1 and Remark 2 in \cite{Martel})

\medskip

\noindent\textbf{Asymptotic $N$-soliton-like solution \cite{Martel}.}\textit{
Let $N\geq 1$ and $x_1,\ldots,x_N\in \mathbb{R}$.
Let $0<c_N<\ldots<c_1<c_*(f)$ be such that \eqref{stable} holds for all $c_j$, $j=1,\ldots,N$.
Then, there exists a unique $H^1$ solution $u(t)$ of \eqref{kdvf} such that}
$$
\lim_{t\to -\infty} \Big\|u(t)-\sum_{j=1}^N Q_{c_j}(.-x_j-c_j t)\Big\|_{H^1}=0.
$$

\smallskip

\noindent Recall also that this behavior is in some sense stable in the energy space,
see Martel, Merle and Tsai \cite{MMT}.

\medskip

We are concerned with the problem of   collision of two solitons. This is a classical 
problem in nonlinear wave propagation which we briefly review (see also the introduction of \cite{MMcol1} and to the references therein).
First, Fermi, Pasta and Ulam \cite{FPU} and Zabusky and Kruskal  \cite{KZ} have exhibited from the numerical point of view
remarkable phenomena related to soliton collision. Next, Lax  \cite{LAX1}  has developed a mathematical framework to study these problems,
known now as complete integrability. The   inverse scattering transform   (for a review on this theory, we refer  for example to
Miura \cite{Miura}) then provided explicit formulas for 
$N$-soliton solutions (Hirota \cite{HIROTA}): let $f(u)=u^2$ or $f(u)=u^3$,
and let $c_1>\ldots>c_N>0$, $\delta_1,\ldots,\delta_N\in \mathbb{R}$. There exists an explicit  solution $U(t,x)$ of \eqref{kdvf} which satisfies
	\[ 	\biggr\|U(t,x)- \sum_{j=1}^N Q_{c_j}(.-c_jt-\delta_j)\biggr\|_{H^1} \mathop{\longrightarrow}_{t\to -\infty} 0,\quad
	\biggl\|U(t,x)- \sum_{j=1}^N Q_{c_j}(.-c_jt-\delta_j')\biggr\|_{H^1}
	\mathop{\longrightarrow}_{t\to +\infty} 0,
\]
for some $\delta_j'$ such that the shifts $\Delta_j=\delta_j'-\delta_j$ depend on the $(c_k)$. For example, the following function $U_{1,c}$, is a $2$-soliton solution of \eqref{kdvf} with $p=2$, ($0<c<1$):
	\begin{equation}\label{NSOL}
	 U_{1,c}(t,x)=6 \frac {
	\partial^2}{
	\partial x^2} \log\bigl(1+e^{x-t}+e^{\sqrt{c}(x-ct)}+\alpha e^{x-t}e^{\sqrt{c}(x-ct)}\bigr) \quad \textrm{with}\quad  \alpha=\left(\frac {1-\sqrt{c}}{1+\sqrt{c}}\right)^2.\end{equation}
As pointed out in \cite{MMcol1}, the problem of describing the collision of two traveling waves is  a general problem for nonlinear PDEs, which is completely open, except in the
integrable case described above.
This kind of problems have been studied since the 60's from both experimental and numerical points of view.

We recall some numerical works for equations of gKdV type.   Bona et al. \cite{BPS}, and
Kalisch and Bona \cite{KB}, studied numerically the problem of collision of two solitary waves for the Benjamin   and the BBM equations.
Shih \cite{SHIH} studied the case of the gKdV equation \eqref{kdvf} with some 
half-integer values of $p$.  Li and Sattinger \cite{LS} investigated the collision problem for the
Ion Acoustic Plasma equation, and Craig et al. \cite{Craig} report on numerics for the Euler equation with free surface. In all these works, the numerics match the experiments and show that for these models, unlike for
the pure solitons of the integrable case, the collision of two solitary waves fails to be elastic
by a very small but non zero dispersion.
 
Finally, 
the multi-soliton solutions
of the NLS (nonlinear Schr\"odinger) model, with special nonlinearity and under spectral assumptions (ruling out the existence of small solitary waves) have been studied by Perelman \cite{P} and
Rodnianski, Schlag and Soffer \cite{RSS}
(in a special case
where the collision has a negligeable effect on the solitary waves due to a very small time of interaction). See also Cao and Malomed \cite{CM}, 
Holmer,  Marzuola and  Zworski \cite{HMZ},
and Holmer and  Zworski \cite{HZ} for the case of the collision of a soliton of the
NLS equation with a Dirac.

\medskip

In \cite{MMcol1}, we present a complete rigorous description of the collision of 
two solitons of \eqref{kdvf} for the nonlinearity $f(u)=u^4$ 
in the case where one soliton is small with respect to the other. First, we prove  that
the collision is not completely elastic in this case i.e. there
does not exist pure $2$-soliton solution (Theorem 1.1 in \cite{MMcol1}).
Note that this is the first rigorous result related to inelastic (but close to elastic) collision,
and that a precise measurement of the defect follows from the analysis (see Theorems 1.1 and 1.2
in \cite{MMcol1}).
 We also
prove that for any solution behaving as $t\to -\infty$ approximately as the sum of two solitons of different sizes,
the two solitons are preserved after the collision,
with a residual term   very small compared to the sizes of the two solitons. Moreover, we give a detailled
description of the collision such that explicit formulas for the main orders of the shifts on the trajectories of the solitons (see Theorems 1.2 and 1.3 in \cite{MMcol1}).

\medskip
 
In this paper, we consider the same questions for \eqref{kdvf} with a general nonlinearity
$f(u)$ satisfying \eqref{surf}.
We consider two solitons 
$Q_{c_1}>0$, $Q_{c_2}$, where the condition on $c_1$, i.e. $0<c_1<c_*(f)$ is not 
restritive,  indeed for many nonlinearities $c_*(f)=+\infty$. In fact, 
in higher dimensions, $c_*(f)=+\infty$  is a typical assumption to study
the   generalizations of \eqref{ellipticf}.
Concerning $Q_{c_2}$, we assume $c_2$ small (depending on $c_1$). In particular, for
$p=3$, we have both a positive and a negative solution. Theorems below apply 
to both solutions.

 Our approach   is the same as in \cite{MMcol1}, the main tool
being the construction of an approximate solution in the collision region. 
The large time behavior is controlled by
asymptotic arguments, from \cite{MM1}, \cite{MMT}, \cite{MMnonlinearity} later refined in 
\cite{yvanSIAM}, \cite{MMas1} and  \cite{MMas2}.
Our first result concerns the asymptotic $2$-soliton like solution at $-\infty$.

\begin{theorem}[Behavior after collision of the asymptotic $2$-soliton-like solution]\label{PURE}~\\
Let $p=2$, $3$ or $4$.
    Assume that $f$ satisfies \eqref{surf}.
    Let $0<c_1<c_*(f)$ be such that the positive solution $Q_{c_1}$ of
    \eqref{ellipticf} with $c=c_1$
    satisfies \eqref{stable}.  
    There exist $c_0=c_0(c_1)\in (0,c_1)$ and $K=K(c_1)>0$ such that  
    for any $0<c_2<c_0$, if $Q_{c_2}$ is a solution of \eqref{ellipticf} with $c=c_2$, then the following holds.
    Let $u(t)$ be the solution of \eqref{kdvf} satisfying
    \begin{equation}\label{th01}
    \lim_{t\to -\infty} \|u(t)-Q_{c_1}(.-c_1 t)-Q_{c_2}(.-c_2 t)\|_{H^1} =0.
    \end{equation}
    Then, there exist  $\rho_1(t),\rho_2(t)$, $c_1^+>c_2^+>0$ and  $K>0$ such that
    \begin{equation*}
    w^+(t,x)=u(t,x)-Q_{c_1^+}(x-\rho_1(t))-Q_{c_2^+}(x-\rho_2(t)) \quad
        \end{equation*}
        satisfies $\sup_{t\in \mathbb{R}} \|w^+(t)\|_{H^1}\leq K c^{\frac 1{p-1}}$ and
       \begin{equation}\label{th02}
        \lim_{t\to +\infty} \|w^+(t)\|_{H^1(x>\frac 1{10} c_2 t)}=0,\quad
        \limsup_{t\to +\infty} \|w^+(t)\|_{H^1} \leq K c^{\frac 2{p-1} - \frac 14 -\frac 1{100}}.
   \end{equation}
   \begin{equation}\label{th02tri}
    \lim_{t\to +\infty}|\rho_1'(t)-c_1^+|+|\rho_2'(t)-c_2^+|=0.
   \end{equation}
   Moreover, $\lim_{t\to +\infty} E(w^+(t))=E^+$ and  $\lim_{t\to +\infty} \int (w^+)^2(t)=M^+$
   exist and 
   \begin{equation}\label{th02bis}
   \frac 12 \limsup_{t\to +\infty}  \int ((w^+_x)^2  +   c_2  (w^+)^2)(t)\leq 2 E^+ + c_2 M^+ \leq \liminf_{t\to +\infty} \int ((w^+_x)^2 + 2 c_2 (w^+)^2)(t),
   \end{equation}
     \begin{equation}\label{th03}
            \frac 1{K} (2 E^+ + c_2 M^+) \leq \frac {c_1^+}{c_1} -1
            \leq {K}   (2 E^+ + c_2 M^+),
     \end{equation}
        \begin{equation}\label{th04} 
        \frac 1K c_2^{\frac 2{p-1}-\frac 12} (2 E^+ + c_1 M^+) 
          \leq 1-\frac {c_2^+}{c_2} 
            \leq {K} c_2^{\frac 2{p-1}-\frac 12} (2 E^+ + c_1 M^+) .
   \end{equation}
\end{theorem}
By time and space translation invariances, the conclusions of Theorem \ref{PURE} hold for any
asymptotic $2-$soliton solution.
If $p=2$ or $4$, $Q_{c_2}$ is necessarily a positive solution. If $p=3$,
$Q_{c_2}$ can be positive or negative.
By considering $-f(-u)$ instead of $f(u)$
one can also consider the case   $Q_{c_1}<0$ for $p=3$.

\medskip

\noindent\textbf{Remark 1}

 1. 
Note that there exists $K>0$ such that for $c$ small, 
\begin{equation}\label{INQc}
\forall x\in \mathbb{R},\quad
\frac 1K c^{\frac 1{p-1}}e^{-\sqrt{c} |x|} \leq |Q_c(x)|\leq K c^{\frac 1{p-1}} e^{-\sqrt{c} |x|},
\end{equation}
so that, for $c$ small,
\begin{equation}\label{ordres}
\|Q_c\|_{H^1}\sim K_1 c^{\frac 1{p-1}-\frac 14},\quad 
\|Q_c\|_{L^\infty}\sim K_2 c^{\frac 1{p-1}}.
\end{equation}
A main information provided by Theorem \ref{PURE} is that the $2$-soliton structure is preserved 
for all time at the main order.
Indeed, we observe that for $p=2,3$ or $4$, $\frac 2{p-1} - \frac 14 -\frac 1{100}<
\frac 1{p-1}$,  thus  from \eqref{th02} and \eqref{ordres},  the two soliton structure is recovered asymptotically in large time. Moreover, 
since $\sup_{t\in \mathbb{R}} \|w^+(t)\|_{H^1}\leq K c_2^{\frac 1{p-1}}\ll \|Q_{c_2}\|_{H^1}$, the $2$ soliton
structure is preserved also during the collision. Note that this estimate is
optimal, the perturbation due to the collision being exactly of size $c_2^{\frac 1{p-1}}$  in $H^1$ during the collision region.

 Theorems \ref{EXISTf} and \ref{STABf} below give   other illustrations of the stability of the two soliton dynamics through the collision.

\smallskip

2. Estimate \eqref{th03} means that the speed of the soliton $Q_{c_1}$ can only increase through the interaction,
and that if $c_1^+=c_1$ then $u(t)$ is   a pure $2$-soliton solution both at $+\infty$ and $-\infty$.  Similarly, $c_2$ can only decrease. Remarkably, for $p=3$,
the property does not depend on the sign of $Q_{c_2}$.

\medskip

Note that it is well-known for the case $f(u)=u^2$ or $u^3$ that the solution $u(t)$ 
considered in Theorem \ref{PURE} is pure at $\pm\infty$
($u(t)$ is explicit in the integrable cases). In contrast, 
in the case $f(u)=u^4$ it was proved in Theorem 1.1 of \cite{MMcol1} that there exists no pure $2$-soliton solution
at both $+\infty$ and  $-\infty$. In the general case $f(u)$, whether
or not the collision is elastic is an open question. A natural question related to
Theorem \ref{PURE} is thus to try to understand, in the case of a general nonlinearity $f(u)$
in which situation the collision is elastic or inelastic, and what is the size of the defect.

\medskip

Our second result is related to the construction of an object similar to the $2$-soliton
solutions with a perturbation term, such that the speeds as $t\to \pm\infty$
are the same. We also obtain an explicit formula for the first order of the resulting shift
on the first soliton. The formula is related to the functions
$c\mapsto \int Q_c$ and $c\mapsto \int Q_c^2$ for $c$ close to $c_1$.

\begin{theorem}[Existence of $2$-soliton like solutions]\label{EXISTf}
Let $p=2$, $3$ or $4$.
    Assume that $f$ satisfies \eqref{surf}.
    Let $0<c_1<c_*(f)$ be such that the positive solution $Q_{c_1}$ of \eqref{ellipticf} with $c=c_1$ satisfies \eqref{stable}.  
    There exist $c_0=c_0(c_1)\in (0,c_1)$ and $K=K(c_1)>0$ such that  
    if $0<c_2<c_0$, and $Q_{c_2}$ is solution of \eqref{ellipticf} with $c=c_2$,
    then there exist a global $H^1$ solution
    $\varphi(t)=\varphi_{c_1,c_2}(t)$ of \eqref{kdvf} and $\Delta_1, \, \Delta_2\in \mathbb{R}$, $\rho_1(t)$, $\rho_2(t)$ satisfying, for all $t,x\in \mathbb{R}$,
    \begin{equation}\label{TH1Z}
   		\varphi(-t,-x)=\varphi(t ,x ),
		\end{equation}
and such that the following holds for $w^{\pm}(t)$ where
\begin{equation*}
w^{-}(t,x)=\varphi(t,x) - Q_{c_1}(x+\rho_1(-t)) - Q_{c_2}(.+\rho_2(-t)),
\end{equation*}
\begin{equation*}
w^{+}(t,x)=\varphi(t,x) - Q_{c_1}(x-\rho_1( t)) - Q_{c_2}(.-\rho_2( t)),
\end{equation*}
    \begin{enumerate}
        \item Asymptotic behavior at $\pm\infty$:
                \begin{equation}\label{TH1A}
          \lim_{t\to -\infty}
            \|w^-(t)\|
            _{H^1(x< \frac{c_2 t}{10})} =0,\quad  \lim_{t\to +\infty}
            \|w^+(t)\|
            _{H^1(x>\frac{c_2 t}{10})} =0,        \end{equation}
            \begin{equation}\label{TH1AAA}
            \lim_{t\to +\infty} |\rho_1'(t)-c_1|+|\rho_2'(t)-c_2|=0.
            \end{equation}
     \item Distance to the sum of two solitons: there exists $t_0>0$ such that
   \begin{equation}\label{TH1E}\begin{split}
        &		 	\|w^-(-t)\|_{H^1} +	  \|w^+(t)\|_{H^1} 
            \leq  K   c_2^{\frac 2{p-1}-\frac 14 -\frac 1{100} }, \quad\text{for all $t>t_0$.}
            \end{split}
\end{equation}
\item Shift property: there exist $\delta_1(c_1)$, $\delta_2(c_1)\in \mathbb{R}$ such that
for $
T_{c_1,c_2}=c_1^{\frac 32} \big(\frac {c_2}{c_1}\big)^{-\frac 12-\frac 1{100}},$
\begin{equation}\label{TH1BBB}
|\rho_1(T_{c_1,c_2})- (c_1 T_{c_1,c_2} +\tfrac 12 \Delta_1)|\leq
K  c_2^{\frac 2{p-1}-\frac 12},\quad
|\rho_2(T_{c_1,c_2})- (c_2 T_{c_1,c_2} +\tfrac 12 \Delta_2)|\leq
Kc_2^{\frac 1{12}}, 
\end{equation}
		\begin{equation}\label{TH1B}
			\left|\Delta_1 - \left(\frac{c_2}{c_1}\right)^{\frac 1{p-1}-\frac 12} \delta_1(c_1)
			  \right|\leq K  c_2^{\frac 2{p-1}-\frac 12},\quad
            \left|\Delta_2- \delta_2(c_1) \right|\leq K c_2^{\frac 1{12}}.
            \end{equation}
            Moreover,
            \begin{equation}
\label{TH1XX}
\delta_1(c_1)= 2\, \mathrm{sgn}(Q_2(0))\,
		 \frac {\int Q_{c_1}\,\frac d{dc} \int {Q_c} _{| c=c_1}} {\frac d{dc} \left(\int {Q_c^2}\right)
			_{| c=c_1}}.
\end{equation}
    \end{enumerate}
\end{theorem}

\noindent\textbf{Remark 2}

1. By \eqref{surf}, assuming $c_1$ small is sufficient to ensure the assumptions of   Theorem~\ref{EXISTf}.
However, Theorem \ref{EXISTf} holds for any
$(c_1,c_2)$ such that $0<c_1<c_*$,  $0<c_2<c_0(c_1)$ and \eqref{stable} holds for $c_1$.

\smallskip

2. 
Recall that   $\|Q_{c_2}\|_{L^2}\sim K c_2^{\frac 1{p-1} -\frac 14}$.
This is to be compared with the size
of $w^\pm(t)$ in \eqref{TH1E}.
Note that in estimate \eqref{TH1E}, $\frac 1{100}$ has no particular meanning.
By the technique of the present paper, one can get   
$\|w^+(t)\|_{H^1}\leq  K(\epsilon_0)   c_2^{\frac 2{p-1}-{\frac 14}-\epsilon_0 }$,
for any $\epsilon_0>0$,  which is sharp, see a lower bound on 
$w^+(t)$ for the case $f(u)=u^4$, in Theorem 1.2 in \cite{MMcol1}.

\smallskip

3. If there exists a Viriel property for $f(u)$ and $Q_{c_1}$, as it is the case for
$f(u)=u^p$ ($p=2,3,4$, see \cite{MMcol1}, \cite{MMas2}), then
$\rho_j(t)-c_j t\to x_j^+$ as $t\to +\infty$, for some $x_j^+$ ($j=1,2$).
In particular, it is the case if $c_1$ is small since then the problem is a pertubation of
$f(u)=u^p$ and the Viriel argument still works for $f(u)$.

Note also that at $t=T_{c_1,c_2}$, the two solitons are already decoupled, by exponential decay.
Thus, \eqref{TH1BBB} means that through the collision, the two solitons are
shifted by $\Delta_1$, respectively, $\Delta_2$ at the first order.
In \eqref{TH1B}, we see that the main part of $\Delta_1$ (if $\delta_1\neq 0$) is the product of
a power of $\frac {c_2}{c_1}$ (depending only on $p$)
by $\delta(c_1)$ which depends on $Q_{c_1}$ and thus
 on the nonlinearity $f(s)$ on the interval $s\in [0,Q_{c_1}(0)]$.
By the stability assumption, we have $\frac d{dc} \int {Q_c^2}_{|_{c=c_1}}>0$,
but the other term in \eqref{TH1XX} $\frac d{dc} \int {Q_c}_{|_{c=c_1}} $ may have 
any sign (for example, for $f(u)=u^p$, $p=2$, $3$ and $4$ this term is respectively 
positive, zero and negative, see \cite{MMcol1}).
Note that the shift  on $Q_{c_1}$ depends on the sign of $Q_{c_2}$.

Similarly, we observe that  $\delta_2(c_1)$  depends only on $c_1$. Thus, if
$\delta_2\not =0$, it follows that the main order of the shift on $Q_{c_2}$ is independent of $c_2$.
In \cite{MMcol1}, we have computed $\delta_2$ for $f(u)=u^4$ and there are well-known formulas
for the case $p=2$, $3$ (see e.g. Miura \cite{Miura}).

\begin{theorem}[Stability of the $2$-soliton structure]\label{STABf}
            Let  $\varphi(t)=\varphi_{c_1,c_2}(t)$ be the solution constructed
    in Theorem~\ref{EXISTf}, under the same assumptions.
    There exists $c_0=c_0(c_1)\in (0,c_1)$ and $K=K(c_1)>0$ such that  
    if $0<c_2<c_0$
    then  the following holds.
 Assume that 
    \begin{equation}\label{hyppp} 
	\| u_0 -    \varphi(0)\|_{H^1}
	\leq     c^{\frac 1{p-1}+\frac 12}_2,
  \end{equation}
  and let   $u(t)$ be the $H^1$ solution of \eqref{kdvf}. 
    Then,  there exist  $  \rho_1(t),\, \rho_2(t) \in \mathbb{R}$ and
    $c_1^\pm,\,c_2^\pm >0$ such that
    \begin{enumerate}
        \item Global in time stability:   $$w(t,x)=u(t,x) - Q_{c_1}(x-  \rho_1(t))
- Q_{c_2}(x- \rho_2(t)) \quad \text{satisfies}$$
        \begin{equation}\label{stabbb}
             \|w(t)\|_{H^1}
\leq K    c_2^{\frac 1{p-1}} ,
\quad \text{for all $t\in \mathbb{R}.$}                     
        \end{equation}
        \item Asymptotic stability: 
        \begin{equation*} 
         \lim_{t\to -\infty}
            \|u(t) - Q_{c_1^-}(. -  \rho_1(t)) - Q_{c_2^-}(.- \rho_2(t))\|
            _{H^1(x {<} \frac{c_2 t}{10})} =0,
        \end{equation*}
        \begin{equation*} 
         \lim_{t\to +\infty}
            \|u(t) - Q_{c_1^+}(. -  \rho_1(t)) - Q_{c_2^+}(.- \rho_2(t))\|
            _{H^1(x > \frac{c_2 t}{10})} =0,
        \end{equation*}
        \begin{equation*}
			\left|\frac {c_1^\pm}{c_1} -1\right| \leq K c_2^{\frac 1{p-1}+\frac 12} ,\quad
			\left|\frac {c_2^\pm}{c_2} -1\right| \leq K     c_2^{\frac 14}   .
		\end{equation*}
		\end{enumerate}
        \end{theorem}
Theorem \ref{STABf} is the analogue of Theorem 1.3 in \cite{MMcol1}. Note that since $\|Q_{c_2}\|_{H^1}\sim K c_2^{\frac 1{p-1}-\frac 14}$, \eqref{stabbb} means that the two solitons (even the smaller one) are preserved through the collision. The loss of a power $\frac 12$ in $c$ between \eqref{hyppp} and \eqref{stabbb} is due
to the difference of sizes of $Q_{c_1}$ and $Q_{c_2}$.

\medskip

The paper is organized as follows.
In Section 2, we construct an approximate solution of \eqref{kdvf} in a large time
region including the collision. This section contains  the main new arguments.
 In Section 3, we recall
preliminary results for the asymptotics of the $2$-soliton structure in large time.
In Section 4, we prove Theorems \ref{PURE}, \ref{EXISTf} and \ref{STABf}.

\section{Construction of an approximate $2$-soliton solution}

For the sake of simplicity, we can first assume by scaling that
$c_*(f)>1$  and  
$$c_1=1\quad \text{and}\quad c_2=c<c_0,$$
where $c_0>0$ is to be chosen small enough.
We denote $Q_1=Q>0$ and we suppose that \eqref{stable} holds for $Q$.
Moreover, in what follows, we assume $Q_{c_2}>0$, the case $Q_{c_2}<0$ (and thus $p=3$)
is treated similarly.
We construct an approximate solution of equation \eqref{kdvf} close to the sum of two soliton solutions related to
$Q$ and $Q_c$ on a large time interval containing the collision time.
(The general case
will follow  by a scaling argument, see Corollary \ref{corAVRIL} in section 2.5.)

\smallskip
Let 
\begin{equation}\label{defTc}
    T_c=c^{-\frac 12 -\frac 1{100}}.
\end{equation}
(The power $\frac 1{100}$ in the definition of $T_c$ above can be replaced by any small number, giving a justification of Remark 2 following Theorem \ref{EXISTf}.)

\begin{proposition}[Construction of an approximate solution of the gKdV eq.]\quad \\
\label{AVRILf}
    There exist $c_0(f)>0$ and $K_0(f)>0$ such that for any $0<c<c_0(f)$, there exists
    a function $v=v_{1,c}$ such that  the following hold.
    \begin{enumerate}
        \item Approximate solution on $[-T_c,T_c]$:  for $j= 0, 1,2$,
        \begin{equation}\label{defofS}
        S(t,x)=\partial_{t} v + \partial_{x} (\partial_{x}^2 v - v + f(v)) \quad \text{satisfies}
        \end{equation} 
        \begin{equation}\label{AV0}
            \forall t\in [-T_c,T_c],\quad   \|\partial_x^{j} S(t)\|_{L^2(\mathbb{R})}
                \leq K_0 c^{ \frac 2{p-1}+\frac 34 }.
        \end{equation}
        \item Closeness to the sum of two solitons for $t=\pm T_c$: there exist $\Delta$, $\Delta_c$ such that
\begin{equation}\label{3-16}\begin{split}
 &   \|v(T_c)-Q(.-\tfrac 12\Delta)-Q_c(.+(1-c)T_c-\tfrac 12\Delta_c) \|_{H^1}\leq K_0 c^{ \frac 2{p-1}+\frac 14 },\\&
    \|v(-T_c)-Q(.+\tfrac 12\Delta)-Q_c(.-(1-c)T_c+\tfrac 12\Delta_c) \|_{H^1}\leq K_0 c^{ \frac 2{p-1}+\frac 14 },
\end{split}\end{equation}
                where
        \begin{equation}\label{defDeltabis}
            \left|\Delta - c^{\frac 1{p-1} -\frac 12} \delta \right|\leq K_0 c^{\frac 2{p-1}-\frac 12},\quad
            \left|\Delta_c - \delta_c\right|\leq K_0 c^{\frac 1{12}},
        \end{equation}
        \begin{equation}\label{Deltatri}
        \delta= 2 \frac {\int Q \, \frac {d}{d\widetilde c} {\int Q_{\widetilde c}}_{|\widetilde c=1}}
        {\frac {d}{d\widetilde c} \left({\int Q_{\widetilde c}^2}\right)_{|\widetilde c=1}}.
        \end{equation}
\item Closeness to the sum of two solitons: for all $t\in [-T_c,T_c]$, there exists
$y(t)$ such that
\begin{equation}\label{OUBLI1}
\| v(t) - Q(.-y(t)) - Q_c(.+(1-c)t)\|_{H^1} \leq K_0 c^{\frac 1{p-1}}.
\end{equation}
    \end{enumerate}
\end{proposition}
To prove Proposition \ref{AVRILf}, we follow the same strategy as in \cite{MMcol1}, Sections 2 and 3.
Here, we recall the  main steps and   only mention the parts which have to be adapted. We refer to \cite{MMcol1} for more details.

\medskip

\noindent\emph{Remark.} It follows from the proof of Proposition \ref{AVRILf} that the constants
$c_0(f)$, $K_0(f)$ depend continuouly on $f \in C^{p+4}$.

\medskip

\noindent\emph{Notation.} \quad
For $k,$ $k'$, $\ell$, $\ell'\in \mathbb{N}$, we denote 
\begin{equation*}
    (k',\ell') \prec (k,\ell) \quad \text{if
$k'<k$ and $\ell'\leq \ell$ or if $k'\leq k$ and $\ell'<\ell$.}
\end{equation*}

We denote by $\mathcal{Y}$ the set of functions $g\in C^\infty(\mathbb{R})$ such that
\begin{equation*}
    \forall j\in \mathbb{N},\ \exists K_j,\, r_j>0,\ \forall x\in \mathbb{R},\quad |g^{(j)}(x)|\leq K_j (1+|x|)^{r_j} e^{-|x|}. 
\end{equation*}
Note that $\mathcal{Y}$ is stable by sum, multiplication and differentiation.

\subsection{Choice of a decomposition for $v$}
We look for $v(t,x)$ with a specific structure as in \cite{MMas1}. Let $k_0\geq 1$,
$\ell_0\geq 0$, and 
$$
\Sigma_0 = \{ (k,\ell), ~ 1\leq k\leq k_0, ~ 0\leq \ell\leq \ell_0\}.
$$
We set
\begin{equation*}
    y_{c}=x+(1-c)t  \quad \text{and} \quad R_c(t,x) = Q_c(y_c),
\end{equation*} 
\begin{equation*}
    y=x-\alpha(y_{c}) \quad \text{and} \quad     R(t,x)=Q(y),
\end{equation*}
where for $(a_{k,\ell})_{(k,\ell)\in \Sigma_0}$,  
\begin{equation}\label{defALPHA} 
    \alpha(s)=\int_{0}^s \beta(s') ds',\quad     
    \beta(s)=\sum_{(k,\ell)\in \Sigma_0} a_{k,\ell} \, c^\ell Q_{c}^k(s).
\end{equation}
The form of $v(t,x)$ is
\begin{equation}\label{defv} 
    v(t,x)=Q(y)+Q_{c}(y_{c})+W(t,x),
\end{equation}
\begin{equation}\label{defW}
    W(t,x)=\sum_{(k,\ell)\in \Sigma_0} 
        c^\ell\left(Q_{c}^k(y_{c}) A_{k,\ell}(y)+(Q_{c}^k)'(y_{c}) B_{k,\ell}(y)\right),
\end{equation}
where $ a_{k,\ell} $, $A_{k,\ell}$, $B_{k,\ell}$ are to be determined.

The motivation in \cite{MMcol1} 
for choosing $W$ of the form \eqref{defW} is the stability of the family of functions
\begin{equation}\label{ily}
    \left\{c^\ell Q_c^k,\ c^\ell (Q_c^k)', \ k\geq 1,\ \ell \geq 0 \right\}
\end{equation}
by  multiplication and differentiation due to the power nonlinearity in the equation (see Lemma 2.1 in \cite{MMcol1}). In the case of equation \eqref{kdvf},
for a general nonlinearity
this structure is preserved up to a lower order term (see Lemma \ref{surQc2}). 
Let 
\begin{equation}\label{2.2bis}
S(t,x)=\partial_{t} v + \partial_{x} (\partial_{x}^2 v - v + v^p).
\end{equation}

\begin{proposition}[Decomposition of $S(t,x)$]\label{SYSTEMEf}
Assume that $f$ is of class $C^{k_0+3}$.  Let
\begin{equation}\label{defLy}
    \mathcal{L} w = -\partial_x^2 w + w -f'(Q) w.
\end{equation}
 Then,
    \begin{align*}
            S(t,x) & = 
        \sum_{(k,\ell)\in \Sigma_0}
        c^\ell Q_c^k(y_c)    \Big[a_{k,\ell} (-3 Q+2 f(Q))'(y)    -(\mathcal{L} A_{k,\ell})'(y)\Big]
        \\& \quad
        + \sum_{(k,\ell)\in \Sigma_0}
        c^\ell (Q_c^k)'(y_c)    
        \Big[a_{k,\ell} (-3 Q'')(y) + \left(3A_{k,\ell}'' +f'(Q) A_{k,\ell}\right)(y)    - (\mathcal{L} B_{k,\ell})'(y)\Big]
        \\& \quad
        + \sum_{(k,\ell)\in \Sigma_0}
        c^\ell\left( Q_c^k(y_c)  F_{k,\ell}(y)   +  (Q_c^k)'(y_c) G_{k,\ell}(y)\right)
        +  \mathcal{E}(t,x)
    \end{align*}
    where 
    $F_{k,\ell} $, $G_{k,\ell} $ and $\mathcal{E}$ satisfy, for any
    $(k,\ell)\in \Sigma_0$,
    \begin{itemize}
        \item[{\rm (i)}] Dependence property of $F_{k,\ell}$ and $G_{k,\ell}$:
        The expressions of $F_{k,\ell}$ and $G_{k,\ell}$ depend only on $(a_{k',\ell'})$, $(A_{k',\ell'})$, $(B_{k',\ell'})$ for $(k',\ell') \prec (k,\ell)$.
        \item[{\rm (ii)}] Parity property of $F_{k,\ell}$ and $G_{k,\ell}$:         Assume that for any $(k',\ell')$ such that $(k',\ell')\prec (k,\ell)$
         $A_{k',\ell'}$ is even and    $B_{k',\ell'}$ is odd, then 
        $F_{k,\ell} $ is odd and  $G_{k,\ell}$ is  even.

 Moreover,
$F_{1,0}=(f'(Q))'$ and $G_{1,0} = f'(Q).$
	\item[{\rm (iii)}] Estimate on $\mathcal{E}$: there exists $\kappa(y)>0$
(depending on $(a_{k,\ell})$ and $(A_{k,\ell})$, $(B_{k,\ell})$) such that
\begin{equation}\label{onE}
\forall j=0,1,2,~
\forall (t,x)\in [-T_c,T_c]\times \mathbb{R},\quad 
	|\partial_x^j \mathcal{E}(t,x)|\leq
\kappa(y) (Q_c^{k_0}(y_c)+c^{\ell_0}) Q_c(y_c).
\end{equation}
    \end{itemize}
\end{proposition}

\noindent\emph{Remark.} 
Estimate \eqref{onE} is only a first rough estimate on the rest term, which can not be used without
further information on $\kappa(y)$. 
In Proposition \ref{VandS},   for   the functions $(A_{k,\ell})$, $(B_{k,\ell})$ to be chosen in this paper, we   estimate precisely the size of $\partial_x^j \mathcal{E}$ in $L^2$.

\medskip

Before proving the above proposition, we recall the following properties of $Q_c$, proved in Appendix A.

\begin{lemma}[Properties of $Q_c$]\label{surQc2}
For $0<c\leq 1$, $\forall k,\widetilde k\in \{1,\ldots,k_0\},$
\begin{align}
&
\frac 1K c^{\frac 1{p-1}} e^{-\sqrt{c}|x|} \leq Q_c(x)\leq
K c^{\frac 1{p-1}} e^{-\sqrt{c}|x|},\quad
|Q_c'(x)|\leq K c^{\frac 1{p-1}+\frac 12} e^{-\sqrt{c}|x|},\label{decay}\\
&(Q_c^{k})'(Q_c^{\widetilde k})' = c k \widetilde k Q_c^{k+\widetilde k} 
+  \sum_{p+1\leq k_1 \leq k_0-k-\widetilde k+2} k\widetilde k \,\sigma_{k_1} Q_c^{k+\widetilde k+k_1-2}
+O(Q_c^{k_0+1}),
\label{taylor0}\\ & 
(Q_c^k)''=c k^2  Q_c^k + \sum_{k+p-1\leq k_1 \leq k_0} \sigma^{k*}_{k_1} Q_c^{k_1}+O(Q_c^{k_0+1}),
\label{taylor1}\\
&(Q_c^k)^{(3)}=ck^2  (Q_c^k)' + \sum_{k+p-1\leq k_1 \leq k_0} \sigma^{k*}_{k_1} (Q_c^{k_1})'
+O(Q_c^{k_0+1}),\label{taylor2}\\
&(Q_c^k)^{(4)}=c^2 k^4 Q_c^k + c \sum_{k+p-1\leq k_1 \leq k_0} \sigma^{k**}_{k_1} Q_c^{k_1}
+ \sum_{k+2p-2\leq k_1 \leq k_0} \sigma^{k***}_{k_1} Q_c^{k_1} +O(Q_c^{k_0+1}),\label{taylor3}
\end{align}
where $\sigma_{k_1}$, $\sigma_{k_1}^{k*}$, $\sigma_{k_1}^{k**}$ and $\sigma_{k_1}^{k***}$
are independent of $c$, and where $O(Q_c^k)$  is a function $\mathcal{E}$
satisfying for $j=0,1,2,$
$|\partial_x^{j} \mathcal{E}(t,x)| \leq K Q_c^k(y_c)$, where $K$ is  independent of $c$.
\end{lemma}

\noindent\emph{Proof of Proposition \ref{SYSTEMEf}.}
Inserting $v=R+R_c+W$ in the expression of $S(t,x)$ in  \eqref{2.2bis}, and
using the equations of $R$ and $R_c$, we obtain the following decomposition
(see also \cite{MMcol1}, Proof of Proposition 2.2)

\begin{equation}\label{decS}
    S(t,x)= \mathbf{I} +\mathbf{II} +\mathbf{III}+{\mathbf{IV}},
\end{equation}
where
$$\mathbf{I}=\partial_t R + \partial_x(\partial_x^2 R -R + f(R)),\quad
\mathbf{II}=\partial_x(f(R+R_c)-f(R)-f(R_c)),$$ 
$$ \mathbf{III}= \partial_t W - \partial_x (\overline{\mathcal{L}} W),\quad \text{
where $\overline{\mathcal{L}} W=-\partial_x^2 w + w - f'(R) w$},$$
$$
{\mathbf{IV}}= \partial_x(f(R+R_c+W) -f(R+R_c) -f'(R) W).$$

\smallskip

\noindent Now, we follow exactly the same steps as in Section 2 of \cite{MMcol1}, replacing
Lemma 2.1 in \cite{MMcol1} by Lemma \ref{surQc2} and using Taylor expansions. 
For example,
by \eqref{surf} 
for $k_0\leq p$, we have the following Taylor expansion of $f$ and $F$:
\begin{equation}
\label{taylor}\begin{split}
f(s)=s^p+f_1(s)=s^p+\sum_{p+1\leq k_1\leq k_0} \frac 1{k_1!} s^{k_1} f_1^{(k_1)}(0) + s^{k_0+1}O(1),\\
F(s)=\frac 1{p+1} s^{p+1} +  \sum_{p+2\leq k_1\leq k_0} \frac 1{k_1!} s^{k_1} f_1^{(k_1-1)}(0) +  s^{k_0+1}O(1).
\end{split}\end{equation}

\medskip

\noindent\emph{Decomposition of $\mathbf{I}$.} As in the proof of Lemma A.1 in \cite{MMcol1}, we claim
\begin{align*}
\mathbf{I}&= \beta(y_c) (-3 Q +2 f(Q))'(y) + \beta'(y_c) (-3 Q'')(y) + c \beta(y_c) Q'(y) + \beta''(y_c) (-Q')(y) 
\\&\quad    + \beta^2(y_c) (3 Q^{(3)})(y) + \beta'(y_c) \beta(y_c) (3 Q'')(y) +\beta^3(y_c) (-Q^{(3)})(y)
\\&= \mathbf{I}_1+\mathbf{I}_2+\mathbf{I}_3+\mathbf{I}_4+\mathbf{I}_5+\mathbf{I}_6+\mathbf{I}_7.
\end{align*}
Using Claim \ref{surDD} (Appendix), we deduce that  $\mathbf{I}$ has the following decomposition:
\begin{align}
    \mathbf{I} & = 
    \sum_{(k,\ell)\in \Sigma_0} c^\ell\left( Q_c^k(y_c)  a_{k,\ell} (-3 Q+2 f(Q))'(y)    + (Q_c^k)'(y_c)  a_{k,\ell} (-3 Q'')(y)\right) \label{IF} \\ &
    \quad + \sum_{(k,\ell)\in \Sigma_0}
    c^\ell\left( Q_c^k(y_c)  F_{k,\ell}^{\mathbf{I}}(y)     +    (Q_c^k)'(y_c) G_{k,\ell}^{\mathbf{I}}(y)\right)+O(Q_c^{k_0+1}),
    \label{IG}
    \end{align}
  where the main terms, i.e. \eqref{IF} are coming from $\mathbf{I}_1$ and $\mathbf{I}_2$
  and $F_{k,\ell}^{\mathbf{I}}$, $G_{k,\ell}^{\mathbf{I}}$ satisfy (i)-(ii) of Proposition \ref{AVRILf}.
  
  \medskip
  
\noindent\emph{Decomposition of $\mathbf{II}$.}
For this term, we use the Taylor decomposition of $f$ both at $0$ and at $R$, i.e.
\begin{align*}
f(R+R_c)-f(R)-f(R_c) &=  \sum_{1\leq k_1 \leq p-1} \frac 1{k_1 !} Q_c^{k_1}(y_c) f^{(k_1)}(Q(y))
\\ &+ \sum_{p\leq k_1 \leq k_0} \frac 1{k_1 !} Q_c^{k_1}(y_c) (f^{(k_1)}(Q(y))-f^{(k_1)}(0))
+O(Q_c^{k_0+1}).
\end{align*} 
Then, by
\begin{equation}\label{claimA2}
\partial_x (g(y))=(1-\beta(y_c)) g'(y),
\end{equation}
applied to $g(y,y_c)=f(Q(y)+Q_c(y_c))-f(Q(y))-f(Q_c(y_c))$,
we obtain:
\begin{equation}
\mathbf{II} = 
        \sum_{(k,\ell)\in \Sigma_0}
         c^\ell \left(Q_c^k(y_c)  F_{k,\ell}^{\mathbf{II}}(y)
        + (Q_c^k)'(y_c)  G_{k,\ell}^{\mathbf{II}}(y)\right) +O(Q_c^{k_0+1}),
    \end{equation}  
where $F_{k,\ell}^{\mathbf{II}}$, $G_{k,\ell}^{\mathbf{II}}$ satisfy (i)-(ii).
Note that $F_{1,0}^{\mathbf{II}}=(f'(Q))'$ and
$G_{1,0}^{\mathbf{II}}=f'(Q)$.

\medskip
  
\noindent\emph{Decomposition of $\mathbf{III}$.} 
Since
$    W(t,x)=\sum_{(k,\ell)\in \Sigma_0} 
        c^\ell\left(Q_{c}^k(y_{c}) A_{k,\ell}(y)+(Q_{c}^k)'(y_{c}) B_{k,\ell}(y)\right),
$
we are reduced to compute $\partial_t w - \partial_x (\overline{\mathcal{L}} w)$
for terms of the type 
$w(t,x)=Q_c^k(y_c) A(y)$ and $w(t,x)=(Q_c^k)'(y_c) B(y)$.
We recall (see Claim A.3 in \cite{MMcol1}), for  $A(x)\in C^3$,
    \begin{align*}
            & \partial_t (Q_c^k(y_c) A(y)) - \partial_x (\overline{\mathcal{L}} (Q_c^k(y_c) A(y)))
            \\ & \quad = Q_c^k(y_c) (-\mathcal{L} A)'(y)+(Q_c^k)'(y_c) (3 A''+f'(Q)A - c A)(y)
            \\ & \quad + Q_c^k(y_c) \beta(y_c) (-3 A'' - f'(Q_c) A + c A)'(y)+ Q_c^k(y_c) \beta'(y_c) (-3 A'')(y)
            \\ & \quad  + Q_c^k(y_c) \beta''(y_c) (-A')(y)+ Q_c^k(y_c) \beta^2(y_c) (3 A^{(3)})(y)
            \\ & \quad   + Q_c^k(y_c) \beta'(y_c) \beta(y_c) (3 A'')(y)+ Q_c^k(y_c) \beta^3(y_c) (-A^{(3)})(y)
                         \\ & \quad + (Q_c^k)'(y_c) \beta(y_c) (-6 A'')(y) + (Q_c^k)'(y_c) \beta'(y_c) (-3 A')(y) 
                + (Q_c^k)'(y_c) \beta^2(y_c) (3 A'')(y) 
            \\ & \quad+ (Q_c^k)''(y_c)(3A')(y) + (Q_c^k)''(y_c)\beta(y_c) (-3A')(y)+(Q_c^k)^{(3)}(y_c) A(y).
    \end{align*}
Note that a similar formula holds for $w(t,x)=(Q_c^k)'(y_c) B(y)$ (see Claim A.4 in \cite{MMcol1}).

Then, from Lemma \ref{surQc2} and the decompositions of $\beta(y_c)$, $\beta''(y_c)$, $\beta^2(y_c)$, $\beta'(y_c)\beta(y_c)$ and $\beta^3(y_c)$ (see Claim \ref{surDD}), we obtain the following decomposition for
$\mathbf{III}$:
\begin{align}
\mathbf{III} & = 
    \sum_{(k,\ell)\in \Sigma_0} c^\ell \left(Q_c^k(y_c) (-\mathcal{L} A_{k,\ell})'(y) +
    (Q_c^k)'(y_c) \big( 3 A_{k,\ell}'' + f'(Q) A_{k,\ell} -(\mathcal{L} B_{k,\ell})' \big)(y)\right)
    \label{IIIA}\\ & \quad
    + \sum_{(k,\ell)\in \Sigma_0} 
            c^\ell\left( Q_c^k(y_c) F_{k,\ell}^{\mathbf{III}}(y)    + 
             (Q_c^k)'(y_c)  G_{k,\ell}^{\mathbf{III}}(y)\right)+  
\mathcal{E}_{\mathbf{III}}(t,x)
\label{IIIB}
\end{align}
where $F_{k,\ell}^{\mathbf{III}}$, $G_{k,\ell}^{\mathbf{III}}$ satisfy (i)-(ii)
and $\mathcal{E}_{\mathbf{III}}(t,x)$ satisfies (iii).
  \medskip
  
\noindent\emph{Decomposition of $\mathbf{IV}$.}
Let $\mathbf{N}=f(R+R_c+W) -f(R+R_c) -f'(R) W$. Using   Taylor
formula and \eqref{claimA2}, we obtain
\begin{equation}\label{taylorN} 
  \mathbf{N}  = 
\sum_{k=2}^{k_0}\frac {1}{k!} ((R_c+W)^{k }-R_c^{k }) f^{(k )}(R)+  \mathcal{E}_{\mathbf{N}}(t,x),
\end{equation}
\begin{equation*}
        {{\mathbf{IV}}}=\sum_{\substack{2\leq k \leq k_0 \\ 0 \leq \ell \leq    \ell_0 }} 
        c^{\ell} \left( Q_c^{k}(y_c) F_{k,\ell}^{{{\mathbf{IV}}}}(y) +  (Q_c^{k})'(y_c) G_{k,\ell}^{{{\mathbf{IV}}}}(y)\right)+\mathcal{E}_{\mathbf{IV}}(t,x),
    \end{equation*}
    where $F_{k,\ell}^{\mathbf{IV}}$ and $G_{k,\ell}^{\mathbf{IV}}$ satisfy (i)-(ii)
 and $\mathcal{E}_{\mathbf{IV}}(t,x)$ satisfies (iii).

\subsection{Resolution of the systems $(\Omega_{k,\ell})$}
Proposition \ref{SYSTEMEf} leads to the following decomposition of $S(t,x)$:
\begin{align*}
          &  S(t,x)   = 
        -\sum_{(k,\ell)\in \Sigma_0}
          c^\ell Q_c^k(y_c)    \Big(  (\mathcal{L} A_{k,\ell})' +a_{k,\ell} (3 Q-2 f(Q))'  -F_{k,\ell} \Big)(y)
        \\& \quad
        - \sum_{(k,\ell)\in \Sigma_0}
        c^\ell (Q_c^k)'(y_c)    
        \Big((\mathcal{L} B_{k,\ell})'  +a_{k,\ell} (3 Q'')  - \left(3A_{k,\ell}'' +f'(Q) A_{k,\ell}\right)    -        G_{k,\ell} \Big) (y)
        +  \mathcal{E}(t,x).
    \end{align*}
Therefore, we want to solve by induction on $(k,\ell)$ the following systems:
\begin{equation*}
    (\Omega_{k,\ell})\quad 
    \left\{
    \begin{array}{l}
         (\mathcal{L} A_{k,\ell})' + a_{k,\ell}(3Q -2f(Q))' = F_{k,\ell}\\
         (\mathcal{L} B_{k,\ell})' + a_{k,\ell} (3Q'') - 3 A_{k,\ell}'' - f'(Q) A_{k,\ell}=G_{k,\ell}.
    \end{array}
    \right.
\end{equation*}
The first step is to  establish a general existence result for the model system:
\begin{equation*}
    (\Omega)\quad 
    \left\{
    \begin{array}{l}
         (\mathcal{L} A)' + a(3Q -2f(Q))' = F\\
         (\mathcal{L} B)' + a (3Q'') - 3 A'' - f'(Q) A=G.
    \end{array}
    \right.
\end{equation*}
We introduce some notation and we recall well-known results concerning the operator $\mathcal{L}$.

\begin{claim}\label{surphi}
    The function $\varphi(x)=-\frac {Q'(x)}{Q(x)}$ is odd and    satisfies:
    \begin{itemize}
        \item[{\rm (i)}] $\lim_{x\rightarrow -\infty} \varphi(x)=-1$; $\lim_{x\rightarrow +\infty} \varphi(x)=1$;
        \item[{\rm (ii)}] $\forall x\in \mathbb{R}$, $|\varphi'(x)|+|\varphi''(x)|+|\varphi^{(3)}(x)|\leq C e^{-|x|}$.
        \item[{\rm (iii)}] $\varphi'\in \mathcal{Y}$,  $(1-\varphi^2) \in \mathcal{Y}$.
         \end{itemize} 
\end{claim}

\noindent\emph{Proof of Claim \ref{surphi}.}
By \eqref{surqc}, we have
$\varphi^2=\frac {Q'^2}{Q^2}=1-\frac{2 F(Q)} {Q^2}$, thus (i) is a consequence of \eqref{surf}.
Next, $\varphi'=\frac 1{Q^2}((Q')^2-Q''Q)=\frac 1{Q^2} (Qf(Q)-2F(Q))$, and 
(ii), (iii) follow from \eqref{surf} and the decay of $Q$.

\begin{lemma}[Properties of $\mathcal{L}$]\label{surL}  The operator $\mathcal{L}$ defined in $L^2(\mathbb{R})$ by \eqref{defLy}    is self-adjoint and satisfies the following properties:
    \begin{itemize}
        \item[{\rm (i)}] There exist a unique $\lambda_0>0$, $  \chi_0 \in H^1(\mathbb{R})$,
  $ \chi_0 >0$ such that
  $\mathcal{L}   \chi_0 =-\lambda_0  \chi_0 $.
                \item[{\rm (ii)}] The kernel of $\mathcal{L}$ is 
        $\{\lambda Q', \lambda \in \mathbb{R}\}$.
        Let  $\Lambda Q =\frac d{dc} {Q_{c}}_{| c=1}$, then
         $\mathcal{L} (\Lambda Q)=-Q $.
        \item[{\rm (iii)}] For all   $h \in L^2(\mathbb{R})$ such that $\int h Q'=0$, 
        there exists a unique $\widetilde h \in H^2(\mathbb{R})$  such that $\int \widetilde hQ'=0$ and $\mathcal{L} \widetilde h=h$; moreover,
        if $h$ is even (respectively, odd), then $\widetilde h$ is even (respectively, odd).
        \item[{\rm (iv)}] For $h\in H^2(\mathbb{R})$,  $\mathcal{L} h \in \mathcal{Y}$ implies $h\in \mathcal{Y}$.
	\item[{\rm (v)}]  If $\frac d{d\widetilde c} {\int Q_{\widetilde c}^2}_{|\widetilde c=c}>0$
then there exists  $\lambda_c>0$ such that
$$
\int  w Q_c=\int w Q_c'=0\quad \Rightarrow \quad
\int (w_x^2+c w^2 -f'(Q_c) w^2) \geq \lambda_c \int w^2.
$$
    \end{itemize}
\end{lemma}

\medskip

\noindent\emph{Proof of Lemma \ref{surL}.} See Weinstein \cite{We1} and proof of Lemma 2.2 in \cite{MMcol1}.

\medskip

We claim the following general existence result for $(\Omega)$ (similar to Proposition 2.3 in \cite{MMcol1}):

\begin{proposition}[Existence for the model problem $(\Omega)$]\label{SYS4}
    Let $F,G:\mathbb{R}\to\mathbb{R}$ such that
    \begin{align*}
        F(x) & = \overline{F}(x) + \widetilde{F}(x) + \varphi(x) \widehat{F}(x),\quad
        G(x)   = \overline{G}(x) + \widetilde{G}(x) + \varphi(x) \widehat{G}(x),
    \end{align*} 
    \begin{itemize}
        \item $\overline{F}$, $\overline{G}\in \mathcal{Y}$;
          $\overline{F}$ is odd and $\overline{G}$ is even;
        \item $\widetilde{F}$ and $\widehat{G}$ are odd polynomial functions;
         $\widehat{F}$ and $\widetilde{G}$ are even polynomial functions.
    \end{itemize}
    Then, there exist $a\in \mathbb{R}$ and two functions $A(x)$, $B(x)$ satisfying
    $(\Omega)$ and such that
    \begin{align*}
        A(x) & = \overline{A}(x) + \widetilde{A}(x) + \varphi(x) \widehat{A}(x), \quad 
        B(x)   = \overline{B}(x) + \widetilde{B}(x) + \varphi(x) \widehat{B}(x),
    \end{align*}
    \begin{itemize}
        \item $\overline{A}$, $\overline{B}\in \mathcal{Y}$;
          $\overline A$ is even and  $\overline B$ is odd;
        \item $\widetilde{A}$ and $\widehat{B}$ are even polynomial functions;
        $\widehat{A}$ and $\widetilde{B}$ are odd polynomial functions.
    \end{itemize}
    Moreover,
    \begin{align}
        &\text{if $\widetilde F=0$ (respectively, $\widehat F=0$) then $\widetilde A=0$ (respectively, $\widehat A=0$);}
                 \label{DEGA1}\\
        &\text{if $\widetilde A''=0$ and $\widetilde G=0$ then $\widetilde B=0$;}\quad 
        \text{if $\widehat A''=0$ and $\widehat G=0$ then $\deg \widehat B=0$.}\label{DEGB2}
    \end{align}
\end{proposition}

\noindent\emph{Remark.} In Proposition \ref{SYS4}, we find one solution of system $(\Omega)$. 
This solution is not unique but this does not play a role in this paper.
See Corollary 3.1 in \cite{MMcol1} for the uniqueness question.

Note that as a consequence of \eqref{DEGB2}, it could be that $\widehat B=b \in \mathbb{R}$
while $\widehat A''= \widehat G=0$. This has the consequence to possibly develop polynomial growths
in the functions $A_{k,\ell}$, $B_{k,\ell}$. In the rest of this paper,  it will
be sufficient to consider indices $(k,\ell)$ for which $\widehat B_{k,\ell}$ is a  constant and the other polynomials
$\widetilde A,\widehat A=0,\widetilde B=0 $ are zero, see Proposition \ref{SYSkl4}. 
However, if one wants to solve the systems $(\Omega_{k,\ell})$ for large $k,\ell$, polynomial growths appear in 
general, see \cite{MMcol1}.

\medskip

\noindent\emph{Sketch of the proof of Proposition \ref{SYS4}.}
As in the proof of Proposition 2.3 in \cite{MMcol1}, we first reduce the proof 
to the case where the second members do not contain polynomials and thus are in $\mathcal{Y}$.

\medskip

\emph{Step 1.} 
Following step 1 of the proof of  Proposition 2.3 in \cite{MMcol1}, considering
\begin{align*}
    &    -\widetilde{A}''(x)+\widetilde{A}(x) = \int_0^x \widetilde{F}(z) dz,
    \quad    -\widehat{A}''(x)+\widehat{A}(x) = \int_0^x \widehat{F}(z) dz,
\end{align*}
\begin{equation*}
  -\widetilde{B}''(x){+}\widetilde{B}(x) = \int_0^x \left(\widetilde{G}(z)  {+} 3 \widetilde{A}''(z)\right)dz,
    \quad   -(\widehat{B}^*)''(x){+}\widehat{B}^*(x) = \int_0^x \left(\widehat{G}(z)  {+} 3 \widehat{A}''(z)\right)dz,
\end{equation*}
where $\widehat{B}=\widehat{B}^*+b$, 
and using  the exponential decay of $f'(Q)$, we reduce ourselves to solving
 the following system in $(a,b,\overline{A},\overline{B})$:
\begin{equation*}
    \left\{
    \begin{array}{ll}
         (\mathcal{L} \overline{A})' + a (3Q -2 f(Q))' = \mathcal{F}                 \\
         (\mathcal{L} \overline{B})' + a (3Q'') - 3 \overline{A}'' - f'(Q) \overline{A}=\mathcal{G} + b (\mathcal{L} \varphi)', 
    \end{array}
    \right.
\end{equation*}
where $\mathcal{F}\in \mathcal{Y}$ is odd,  $\mathcal{G}\in \mathcal{Y}$ is even
and $\mathcal{F}$, $\mathcal{G}$ do not depend on the  parameters $a$ and $b$.
See \cite{MMcol1} for more details.

\medskip

\emph{Step 2.} Existence of a solution to the reduced  system. 
Set $   \mathcal{H}(x)=\int_{-\infty}^x \mathcal{F}(z) dz.$
Since $\mathcal{F}$ is odd, $\int_{\mathbb{R}} \mathcal{F}=0$ and so
$\mathcal{H}\in \mathcal{Y}$ is even.
To find a solution $(a,b,\overline A,\overline B)$ of $(\overline{\Omega})$, it is sufficient to solve
\begin{equation*}
    (\overline{\Omega})\quad 
    \left\{
    \begin{array}{ll}
         \mathcal{L} \overline{A} + a (3Q -2 f(Q))= \mathcal{H}  \\
         (\mathcal{L} \overline{B})' + a (3Q'') - 3 \overline A'' -f'(Q) \overline A=\mathcal{G} + b (\mathcal{L} \varphi)'.
    \end{array}
    \right.
\end{equation*}
Since $\int \mathcal{H}Q'=0$ (by parity) and $\mathcal{H}\in \mathcal{Y}$, it follows from Lemma \ref{surL} (iii)-(iv) that there exists  
\begin{equation}\label{defofH}
    \overline{H} \in \mathcal{Y} , \text{ even, such that } \mathcal{L}\overline{H}=\mathcal{H}.
\end{equation}
By Lemma \ref{surL}, there also exists
\begin{equation}\label{defofV0}
     V_0 \in \mathcal{Y}, \text{ even, such that } \mathcal{L} V_0=3 Q -2 f(Q).
\end{equation}
It follows that, for all $a$,
\begin{equation}
    \overline{A} = \overline{H} - a V_0
\end{equation}
is solution of $\mathcal{L} \overline{A} + a(3 Q - 2 f(Q))= \mathcal{H}$, 
moreover, $\overline{A}$ is even and $\overline A\in \mathcal{Y}$.  
Note that at this point $(a,b)$ are still free, they will be used to  solve the second equation. 
Indeed,  replacing $\overline{A}$ by $\overline{H} - a V_0$ in this equation, solving
$(\overline{\Omega})$ is equivalent  to finding $(a,b,\overline{B})$ such that
\begin{equation}\label{interB} 
    (\mathcal{L} \overline{B})' =- a Z_0 + D + b (\mathcal{L} \varphi)',
\end{equation}
where
\begin{equation*}
    D= 3 \overline{H}''+ f'(Q) \overline{H}+\mathcal{G},
\quad
    Z_0=3 Q'' + 3 V_0'' + f'(Q) V_0.
\end{equation*}
It follows from the properties of $Q$, $V_0$, $\mathcal{G}$ and $\overline{H}$ that $D$ and $Z_0$ are even and satisfy
$Z_0$, $D\in \mathcal{Y}$.
To solve \eqref{interB}, it suffices to find $\overline{B}\in \mathcal{Y}$ such that
\begin{equation}\label{fm}
    \mathcal{L} \overline{B} = E \quad \text{where} \quad E=\int_{0}^x (D-a Z_0)(z) dz + b \mathcal{L}\varphi.
\end{equation}
We now choose $(a,b)$ such that the function $E$ is orthogonal to $Q'$ and has decay at $\infty$.
First, we claim a nondegeneracy condition on $Z_0$, related to the strict stability of the soliton $Q$
(i.e. assumption \eqref{stable}). This is a nontrivial extension of Claim 2.3 in \cite{MMcol1},
which means that the solvability of $(\Omega)$ is related to the noncriticality of $Q$.

\begin{claim}[Nondegeneracy condition]\label{surE} 
            \begin{equation}\label{NONDEGE}
\int Z_0 Q=-\frac 12 \frac d {dc} \int {Q_c^2}_{\big |c=1}=-\int \Lambda Q \, Q\neq 0. 
\end{equation}
\end{claim}

Assuming Claim \ref{surE}, we finish the proof of Proposition \ref{SYS4}.
Let
        \begin{equation}\label{parametres}
            a= \frac {\int DQ}{\int Z_0 Q} \quad \text{and} \quad b= - \int_0^{+\infty} (D-aZ_0)(z) dz.
        \end{equation}
        Then, $E$ defined by \eqref{fm} satisfies
        \begin{equation}\label{decE}
            E\in \mathcal{Y}, \quad E \text{ is odd,} \quad \int E Q'=0 .
        \end{equation}
Indeed, by integration by parts, and decay properties of $Q$, we have
\begin{equation*}
    \int E Q'=-\int (D-a Z_0)Q+ b \int (\mathcal{L} \varphi) Q'= -\int DQ + a \int Z_0 Q
    +b \int \varphi (\mathcal{L} Q')=0,
\end{equation*}
by \eqref{parametres} and $\mathcal{L}Q'=0$.
By Claim \ref{surphi} and \eqref{parametres}, we have
\begin{equation*}
    \lim_{+\infty}E=\int_{0}^{+\infty} \left(D-a Z_0 \right) dz+ b \lim_{+\infty} (\mathcal{L}\varphi)=0
\quad \text{and so $E\in \mathcal{Y}$.}
\end{equation*}
For $(a,b)$ fixed as in    \eqref{parametres}, from \eqref{decE} and Lemma \ref{surL},
it follows that there exists $\overline{B}\in \mathcal{Y}$ such that $\mathcal{L} \overline{B} = E .$
Setting
\begin{equation*}
    A=\overline A+\widetilde A+\widehat A, \quad  B=\overline B+\widetilde B+\widehat B,
\end{equation*}
we have constructed a solution of system $(\Omega)$. \quad \hfill $\Box$

\medskip

\noindent\emph{Proof of Claim \ref{surE}.}  
Let $\Lambda Q$ be defined in Lemma \ref{surL}; recall that $\mathcal{L} (\Lambda Q)= -Q$.
Note also that $\mathcal{L}(xQ')=-2 Q''$ (since $\mathcal{L}Q'=0$).
Thus, $V_0$ defined by \eqref{defofV0} is $V_0=-\Lambda Q-xQ'$. 
Therefore,
\begin{equation*}\begin{split}
\int Z_0 Q & = 3 \int Q'' Q + \int (3 v_0'' + f'(Q) V_0) Q= 
=3 \int Q'' Q + \int V_0(3 Q''+Q f'(Q))
\\& = - 3 \int (Q')^2 - \int (\Lambda Q+xQ') (3Q''+Qf'(Q)).
\end{split}\end{equation*}
First,
\begin{equation*}\begin{split}
-\int xQ' (3 Q''+ Qf'(Q)) & = -\int xQ'(4 Q'' -Q+ f(Q) + Qf'(Q))\\
&=2 \int (Q')^2 - \frac 12 \int Q^2 + \int Q f(Q).
\end{split}\end{equation*}
Since $\mathcal{L} Q= - Q'' + Q -Qf'(Q)$, we also
have $\mathcal{L}(Q+\Lambda Q+xQ')= - 3 Q''- Q f'(Q)$ and thus
\begin{equation*}\begin{split}
-\int \Lambda Q (3 Q''+ Qf'(Q))  &= \int \Lambda Q \mathcal{L} (Q+\Lambda Q+xQ') =- \int Q(Q+\Lambda Q+xQ')\\
& =-\frac 12 \int Q^2 - \int \Lambda Q\, Q.
\end{split}\end{equation*}
Thus, we obtain by $\int (Q')^2 + \int Q^2 =\int Q f(Q)$,
$$\int Z_0 Q = - \int (Q')^2 - \int Q^2 + \int Q f(Q)
-\int \Lambda Q \,Q =- \int \Lambda Q\, Q.$$

\medskip
  Proposition \ref{SYS4}   allows us to solve the systems $(\Omega_{k,\ell})$
  for all $(k,\ell)\in \Sigma_0$,   for any $k_0\geq 1$, $\ell_0\geq 0$
  (as in \cite{MMcol1}).
  In the present paper, for the sake of simplicity, we work for the minimal set of indices so that
  we are able to prove Theorems 1 and 2. Indeed, let us define
  \begin{equation}\label{defSIGMA}
  	\Sigma_p=\{(k,\ell)~ | ~ \ell=0, ~ 1\leq k \leq p, ~ \text{or} ~ \ell=1, ~ k=1\}.
  \end{equation}
 Using Propositions \ref{SYSTEMEf} and \ref{SYS4}, we solve the systems $(\Omega_{k,\ell})$ 
by induction on $(k,\ell)\in \Sigma_p$,
following   \cite{MMcol1}.

\begin{proposition}[Resolution of $(\Omega_{k,\ell})$  for $(k,\ell)\in \Sigma_p$]\label{SYSkl4}
    For all $(k,\ell)\in \Sigma_p$, there exists $(a_{k,\ell},A_{k,\ell},B_{k,\ell})$ of the form
    \begin{equation}\label{ABkl1}
    \begin{split}
        & A_{k,\ell}(x)=\overline A_{k,\ell}(x)\in \mathcal{Y},\quad    B_{k,\ell}(x) = \overline{B}_{k,\ell}(x) + \varphi(x) b_{k,\ell}(x), \quad b_{k,0}\in \mathbb{R},~
           \overline{B}_{k,\ell}\in \mathcal{Y}, \\
           &\text{
          $ {A}_{k,\ell}$ is even and  $ {B}_{k,\ell}$ is odd,}
    \end{split}\end{equation}
    satisfying
    \begin{equation*}
        (\Omega_{k,\ell})\quad 
        \left\{
        \begin{array}{ll}
             (\mathcal{L} A_{k,\ell})' + a_{k,\ell} (3Q -2 f(Q))' = F_{k,\ell}                             &  \\
             (\mathcal{L} B_{k,\ell})' + a_{k,\ell} (3Q'') - 3 A_{k,\ell}'' - f'(Q) A_{k,\ell}=G_{k,\ell}, \qquad         &    
        \end{array}
        \right.
    \end{equation*}
    where  $F_{k,\ell}$, $G_{k,\ell}$ are defined in Proposition \ref{SYSTEMEf}.
\end{proposition}
As a consequence of Proposition \ref{SYSkl4}, we see that by restricting the sum defining $v(t,x)$
to the set of indices
$\Sigma_p$, all the functions $A_{k,\ell}$ belong to $\mathcal Y$ and the functions
$B_{k,\ell}$ are bounded with derivatives in $\mathcal{Y}$. This will   simplify
the proof of the estimates in Proposition \ref{VandS} with respect to the general
estimates proved in \cite{MMcol1}.

\medskip 

\noindent\emph{Proof of Proposition \ref{SYSkl4}.} 
\emph{1. Case $k=1,\ \ell=0$.} Recall that from Proposition \ref{SYSTEMEf},
the functions $F_{1,0},G_{1,0}\in \mathcal{Y}$ are explicit.
Thus, from Proposition \ref{SYS4} \eqref{DEGA1}-\eqref{DEGB2}, the system $(\Omega_{1,0})$ has a solution
$(a_{1,0},A_{1,0},B_{1,0})$ such that
$$
\widetilde A_{1,0}=\widehat A_{1,0}=\widetilde B_{1,0}=0 \quad \text{and} \quad
\widehat B_{1,0}=b_{1,0}, ~ b_{1,0}\in \mathbb{R}.
$$

\medskip

\emph{2. Case $2\leq k\leq p,\ \ell=0$.} In this case, by induction on $1\leq k\leq p$, 
we solve $(\Omega_{k,0})$, and we prove
\begin{equation}\label{recurr}
\widetilde A_{k,0}=\widehat A_{k,0}=\widetilde B_{k,0}=0 \quad \text{and} \quad
\widehat B_{k,0}=b_{k,0}, ~ b_{k,0}\in \mathbb{R}.
\end{equation}
The argument consists in proving that if  property \eqref{recurr} is satisfied for all $1\leq k'<k$,
then $F_{k,0},G_{k,0}\in \mathcal{Y}$, and thus by Proposition \ref{SYS4},  \eqref{recurr} holds for $k$ as
well. This has been checked in   detail in \cite{MMcol1}, see Claim 2.4 and Lemma B1 (except for the case $k=p$).
First, it is quite clear that $\mathbf{I}$ and $\mathbf{II}$ (see Proposition \ref{SYSTEMEf})
contribute of terms 
$F_{k,0}^{\mathbf{I},\mathbf{II}},G_{k,0}^{\mathbf{I},\mathbf{II}}\in \mathcal{Y}$, see also proof of Lemma B.1
in \cite{MMcol1}.
For the term $\mathbf{III}$ in the decomposition of $S(t,x)$, which is linear in $W$, the proof is exactly the same as in Claim 2.4 of \cite{MMcol1}. 

Now, we give some details concerning the term $\mathbf{IV}$. 
Recall first that 
$\mathbf{IV}=\partial_x \mathbf{N}$, where
$\mathbf{N}=f(R+R_c+W)-f(R+R_c)-f'(R)W$.
In the Taylor expansion \eqref{taylorN}, 
for  $2\leq k_1\leq p-1$, the term $f^{(k-1)}(R(x))$ decays as $e^{-|x|}$, by \eqref{surf},
thus the contribution of these terms to $F_{k',\ell'}$, $G_{k',\ell'}$ are in $\mathcal{Y}$.
For $k=p$, the term $f^{(k-1)}(R(x))$ is bounded and the term of 
lower order in $((R_c+W)^p-R_c^p)$ which is not in $\mathcal{Y}$
comes from $B_{1,0}=\overline B_{1,0}+b_{1,0}\varphi$.
Thus, the lowest order term not localized in the $y$ variable is 
$$p b_{1,0}(Q_c^{p-1} Q_c'\varphi)_x
=pb_{1,0} (Q_c^{p-1}Q_c''+(p-1) Q_c^{p-2} (Q_c')^2)
+p b_{1,0}Q_c^{p-1} Q_c'\varphi'.$$ 
Using Lemma \ref{surQc2}, this term does not give contribution
for $\ell=0$, $k=p$.

It follows that
$F_{k,0},G_{k,0}\in \mathcal{Y}$, and thus by Proposition \ref{SYS4}, we obtain a solution
satisfying \eqref{recurr}.

\medskip

\emph{3. Case $k=1$, $\ell=1$.} This case is handled in the same way, we notice
that $F_{1,1},G_{1,1}\in \mathcal{Y}$, and conclude that 
\begin{equation} 
\widetilde A_{1,1}=\widehat A_{1,1}=\widetilde B_{1,1}=0 \quad \text{and} \quad
\widehat B_{1,1}=b_{1,1}, ~ b_{1,1}\in \mathbb{R}.
\end{equation}

\subsection{Definition of $v(t)$ and estimates on $S(t,x)$}

We  define the function $v(t,x)$ as follows.
For $(k,\ell)\in \Sigma_p$, we consider $(a_{k,\ell},A_{k,\ell},B_{k,\ell})$ defined in Proposition
\ref{SYSkl4}, and  $v(t,x)$ defined by
\begin{equation}\label{defvbis} 
    v(t,x)=Q(y)+Q_{c}(y_{c})+ \sum_{(k,\ell)\in \Sigma_p } 
        c^\ell\left(Q_{c}^k(y_{c}) A_{k,\ell}(y)+(Q_{c}^k)'(y_{c}) B_{k,\ell}(y)\right)
        \end{equation}
where $y_c=x+(1-c)t$,
$ y=x-\alpha(y_{c})$ 
and
\begin{equation}\label{defALPHAbis} 
    \alpha(s)=\int_{0}^s \beta(s') ds',\quad     
    \beta(s)=\sum_{(k,\ell)\in \Sigma_p} a_{k,\ell} \, c^\ell Q_{c}^k(s).
\end{equation}
For this choice of function $v(t,x)$ and for $S(t,x)$ defined by \eqref{defofS}, we claim the following estimates.

\begin{proposition}[Estimates on $V$ and $S$]\label{VandS}
    For any  $0<c<1$, 
    for any $t\in   [-T_c,T_c]$, $W(t)$, $S(t)$ belong to $H^s(\mathbb{R})$ for all $s\geq 1$ and satisfy
    \begin{equation}\label{estV}
        \|W(t)\|_{H^1}=\|v(t)-R(t)-R_c(t)\|_{H^1}\leq K c^{\frac 1{p-1}},
    \end{equation}
    \begin{equation}\label{OUBLI3}
    \inf_{y_1\in \mathbb{R}} \|v(t)-Q(.-y_1)-Q_c(.+(1-c)t)\|_{H^1}\leq K c^{\frac 1{p-1}}.
    \end{equation}
    \begin{equation}\label{estS} j=0,1,2,\quad
        \|\partial_x^{(j)}S(t)\|_{L^2}\leq K_j c^{\frac 2{p-1} + \frac 34},
    \end{equation}
    \end{proposition}

\noindent\emph{Proof of Proposition \ref{VandS}.}
The proof of Proposition \ref{VandS} is based on explicit estimates on $|\alpha'|$ and on  all terms of $v(t,x)$ and $S(t,x)$.
Recall from Proposition \ref{SYSkl4} that since $v(t,x)$ is defined only with $(k,\ell)\in \Sigma_p$, we have
$A_{k,\ell}\in \mathcal{Y}$ and $B_{k,\ell}\in L^\infty$, with derivatives in $\mathcal{Y}$.
 
 First, we claim
\begin{equation}\label{estALPHA}
\forall s\in \mathbb{R},\quad 
|\alpha(s)|\leq K c^{ \frac 1{p-1}-\frac 12 },\quad |\beta(s)|=|\alpha'(s)|\leq K c^{\frac 1{p-1}}.
\end{equation}
Indeed, for $c$ small,
\begin{equation*}
    |\alpha(s)|\leq\sum_{ (k,\ell)\in \Sigma_p}\Biggl| a_{k,\ell} \, c^\ell \int_0^s Q_{c}^k(s') ds'\Biggr|\leq
    \max_{(k,\ell)\in \Sigma_p } |a_{k,\ell}| \ 
    \times \sum_{(k,\ell)\in \Sigma_p}     \int Q_{c}^k \leq K   \int Q_{c} .
\end{equation*}
Since $Q_c(s')\leq K c^{\frac 1{p-1}} \exp(-\sqrt{c} |s'|)$, 
$ \|\alpha\|_{L^\infty}\leq K   \int  Q_{c} \leq K c^{\frac{1}{p-1}-\frac 12}. 
$ Similarly,
$    \|\alpha'\|_{L^\infty}\leq K c^{\frac {1}{p-1}}$.

\medskip

Proof of \eqref{estV}. For all $(k,\ell)\in \Sigma_p$, since $A_{k,\ell}\in \mathcal{Y}$ and $B_{k,\ell}\in L^\infty$, we have
\begin{equation*}\begin{split}
&\|c^{\ell} Q_c^k(y_c) A_{k,\ell}(y)\|_{L^2} \leq K c^{\ell} \|Q_c^k\|_{L^\infty} \leq K c^{\frac 1{p-1}},\\
&\|c^{\ell} (Q_c^k)'(y_c) B_{k,\ell}(y)\|_{L^2} \leq K c^{\ell} \|(Q_c^k)'\|_{L^2} \leq K c^{\frac 1{p-1}+\frac 14}.
\end{split}\end{equation*}
The same is true for $\partial_x W(t,x)$ using \eqref{estALPHA}.

\medskip

Proof of \eqref{OUBLI3}. Since $R_c(t)=Q_c(x+(1-c)t)$, we only have to prove that,
for all $t\in [-T_c,T_c]$,
\begin{equation}\label{OUBLI4}
\inf_{y\in \mathbb{R}} \|R(t)- Q(.-y)\|_{H^1}\leq K c^{\frac 1{p-1}}.
\end{equation}
By \eqref{estALPHA}, taking $c$ small enough so that $|\alpha'(t)|<\frac 12$,
 for all $t  \in [-T_c,T_c]$, there exists a unique $y(t)$  such that
$y(t)-\alpha(y(t)+(1-c)t)=0$.
Then, 
\begin{equation*}\begin{split}
\| R(t)-Q(.-y(t))\|_{H^1}
& =\| Q(.-(\alpha(x+y(t)+(1-c)t) - y(t)))-Q\|_{H^1}\\
& = \| Q(.-(\alpha(x+y(t)+(1-c)t) -\alpha(y(t)+(1-c)t)  ))-Q\|_{H^1}
\end{split}\end{equation*}
By \eqref{estALPHA}, we have  $|\alpha(x+y(t)+(1-c)t) -\alpha(y(t)+(1-c)t) |\leq 
K c^{\frac 1{p-1}} |x|$. Thus, we obtain \eqref{OUBLI4}.

\medskip

Proof of \eqref{estS}.
By the decomposition of $S(t,x)$ in the proof of Proposition \ref{SYSTEMEf},
and the choice of $A_{k,\ell}$, $B_{k,\ell}$ in Proposition \ref{SYSkl4},
we obtain $S(t,x)=\mathcal{E}(t,x)$ as defined in Proposition~\ref{SYSTEMEf}.
Thus, we only have to estimate $\mathcal{E}(t)$.
Since for any $(k,\ell)\in \Sigma_p$, $A_{k,\ell}, ~B_{k,\ell}\in L^\infty$ 
(with derivatives in $\mathcal{Y}$), it follows from the decomposition of $S(t,x)$
(see proof of Proposition \ref{SYSTEMEf}) that all functions of the $y$ variable in 
the expression of $S(t,x)$ are bounded. Thus, we have
$$|S(t,x)|\leq K (|Q_c^{p+1}(y_c)|+ c |Q_c^2(y_c)|),$$
where $K>0$ is independent of $y$ and $c$.
Since $\|Q_c^{p+1}(y_c)\|_{L^2}+c \|Q_c^2(y_c)\|_{L^2}\leq K c^{\frac 2{p-1} + \frac 34},$
we obtain
$$
\| \mathcal{E}(t)\|_{L^2}\leq K c^{\frac 2{p-1} + \frac 34}.
$$
The estimates on the derivatives of $S$ are obtained in the same way.

\subsection{Proof of Proposition \ref{AVRILf}}

In what follows, we will see that the first order of the shift $\Delta$ on $Q$ is $a_{1,0} \int Q_c$.
We first derive an explicit formula for $a_{1,0}$ in order to prove Proposition \ref{AVRILf}.

\begin{lemma}[Computation of the first order of the shift on $Q$]\label{COMPUTf}
$$
a_{1,0}=2 \frac { \frac d{d\widetilde c} \int {Q_{\widetilde c}}_{  | \widetilde c=1}}
{\frac d {d\widetilde c} \left(\int {Q_{\widetilde c}^2}\right)_{  | \widetilde c=1}}.
$$
\end{lemma}

\noindent\emph{Proof of Lemma \ref{COMPUTf}.}
From Proposition \ref{SYSTEMEf} and Proposition \ref{SYSkl4}, the system $(\Omega_{1,0})$ writes, for $p=2$, $3$ and $4$:
\begin{equation*}
    (\Omega_{1,0})\quad  \left\{
    \begin{array}{l}
         \mathcal{L} A_{1,0} + a_{1,0} (3Q -2 f(Q)) = f'(Q) \\
         (\mathcal{L} B_{1,0})' + a_{1,0} (3Q'') - 3 A_{1,0}'' - f'(Q) A_{1,0}=f'(Q).
    \end{array}
    \right.
\end{equation*}
Recall from Claim \ref{surE} that
$V_0=-\Lambda Q- x Q'$ solves $\mathcal{L} V_0=3Q - 2 f(Q)$.
Let $V_1$ be the even $H^1$ solution of $\mathcal{L}V_1=f'(Q)$.
Then, the function $A_{1,0}=V_1- a_{1,0} V_0$
solves the first line of $(\Omega_{1,0})$, independently of the value of $a_{1,0}$.
By replacing $A_{1,0}$ in the second line of the system $(\Omega_{1,0})$, we obtain
\begin{equation*}
    (\mathcal{L} B_{1,0})'+ a_{1,0} Z_0=Z_{1},
\end{equation*}
where
\begin{equation}\label{defZs}
    Z_0=3 Q''+ 3 V_0''+f'(Q) V_0, \quad 
    Z_1=3 V_1'' +p Q^{p-1} V_1 + f'(Q).
\end{equation}

Since $\mathcal{L} Q'=0$, we have $\int (\mathcal L B_{1,0})' Q=0$ and so
$$a_{1,0} \int Z_0 Q= \int Z_1 Q.$$
In Claim \ref{surE}, we have obtained
$$\int Z_0 Q=-\int \Lambda Q\, Q=-\frac 12 {\frac d {d\widetilde c} \int {Q_{\widetilde c}^2}_{  | \widetilde c=1}}.$$
Now, we compute $\int Z_1 Q$ similarly as in Claim \ref{surE},
\begin{equation*}
\begin{split}
\int Z_1 Q &= \int Q (3 V_1'' + f'(Q) V_1 + f'(Q))
= \int V_1 (3 Q''+ Qf'(Q)) + \int Q f'(Q)\\
&= - \int \mathcal{L} V_1 (Q+\Lambda Q+xQ') + \int Q f'(Q)
=- \int f'(Q) \Lambda Q + \int f(Q).
\end{split}
\end{equation*}
Now, since $\mathcal{L} (\Lambda Q)=-Q$, we have
$\int \Lambda Q = -\int Q + \int \Lambda Q \, f'(Q)$ and since
$-Q''+ Q = f(Q)$, we have $\int Q= \int f(Q)$.
Thus,
$\int Z_1 Q =-\int \Lambda Q =-\frac d{dc} \int {Q_c}_{| c=1}$, which completes the proof.

\medskip

\noindent\emph{Proof of Proposition \ref{AVRILf}.}
From what precedes (in particular Proposition \ref{VandS}),
we only need to recompose the function $v(t,x)$ at time $\pm T_c$,
combining the first terms of the decomposition of $v(t,x)$.
By symmetry, we consider only $t=T_c$.
This proof follows closely the proof of Proposition 3.1 in \cite{MMcol1}.

\smallskip

\noindent 1. First, we claim  
\begin{equation}\label{recompV}\begin{split}
& \|v(T_c)-Q(y)-Q_c(y_c)-b_{1,0} Q_c'(y_c)\|_{H^1}\leq K c^{\frac 2{p-1} + \frac 14}.
\end{split}\end{equation}
Indeed, from the definition of $v(t,x)$, 
and the fact for $(k,\ell)\in \Sigma_p$, $A_{k,\ell}\in \mathcal{Y}$, $B_{k,\ell}\in L^\infty,$ we have:
$$
|v(T_c)-Q(y)-Q_c(y_c)-b_{1,0} Q_c'(y_c)|\leq
K \left[ Q_c(y_c) e^{-\frac {|y|}2} +
|(Q_c^2)'(y_c)|+c|Q_c'(y_c)|)\right].
$$
By \eqref{decay}, for all $t\in [-T_c,T_c]$,
$
\| Q_c(y_c) e^{-\frac {|y|}2}\|_{H^1} \leq K   \exp(-\tfrac 12 \sqrt{c} t),
$
and thus at $t=T_c$, for $c$ small enough,
$$
\|  Q_c(y_c) e^{-\frac {|y|}2}\|_{H^1}\leq 
K \exp(-\tfrac 12 c^{-\frac 1{100}})\leq K c^{10}.
$$ 
By \eqref{decay}, $\|(Q_c^2)'(y_c)\|_{H^1}+c\|Q_c'(y_c)\|_{H^1}\leq K c^{\frac 2{p-1} + \frac 14}$,
and thus the   estimate is proved for the $L^2$ norm. We proceed similarly
for the estimate on $\partial_x(v(T_c)-Q(y)-Q_c(y_c)-b_{1,0} Q_c'(y_c))$.

\medskip

\noindent 2.  Position of the soliton $Q$ at $t=T_c$.
Let
$$\Delta=\sum_{(k,\ell)\in \Sigma_p} a_{k,\ell} c^\ell \int Q_c^k.$$
We claim
    \begin{align}
    & \text{for $x\geq -T_c/2$ and $t=T_c$},\quad 
                   |\alpha(y_c)-\tfrac 12 \Delta|\leq K    e^{-\frac 14 c^{-\frac 1{100}}},\label{posQ1}\\
    & \text{for $t=T_c$},\quad 
            \|Q(y)-Q(.-\tfrac 12 \Delta)\|_{H^1}\leq K   e^{-\frac 12 c^{-\frac 1{100}}}. \label{posQ2}
        \end{align}

 Proof of \eqref{posQ1}. For any $k\geq 1$, for any $y_c>0$, we have, by \eqref{decay},
\begin{equation*}
    0\leq \int_{y_c}^{\infty} Q_{c}^k(s) ds
    \leq K c^{\frac 1{p-1}}\int_{y_c}^{\infty} e^{-\sqrt{c}\, s } ds=K c^{\frac 1{p-1}-\frac 12} e^{-\sqrt{c}\, y_c},
\end{equation*}
we obtain
\begin{equation*}
    \left|\alpha(y_c)-\tfrac 12 \Delta\right| \leq K c^{\frac 1{p-1}-\frac 12} e^{-\sqrt{c}\,y_c}.
\end{equation*}
For $x\geq -T_c/2$ and $t=T_c$, we have $y_c=x+(1-c)T_c$$\geq (\frac 12 -c)T_c$, thus $\sqrt{c}\, y_c \geq \frac 12 c^{-\frac 1{100}}-1$, and so
\begin{equation*}
    |\alpha(y_c)-\tfrac 12 \Delta|\leq K  c^{-1/6} e^{-\frac 12 c^{-\frac 1{100}}}\leq K   e^{-\frac 14 c^{-\frac 1{100}}} . 
\end{equation*}

Proof of \eqref{posQ2}. 
For $x\geq -T_c/2$, by \eqref{posQ1}, we have $|\alpha(y_c)-\tfrac 12 \Delta|\leq K  c^{ \frac 1{p-1}-\frac 12} e^{-\frac 12 c^{-\frac 1{100}}}$, and so
\begin{equation*}
    \|Q(y)-Q(.-\tfrac 12 \Delta)\|_{H^1(x>-T_c/2)}\leq  K  c  e^{-\frac 14 c^{-\frac 1{100}}}.
\end{equation*}

For $x<-T_c/2$, since $y=x-\alpha(y_c)$, and $|\alpha(y_c)|\leq K c^{\frac 1{p-1}-\frac 12}$, we have $y<-T_c/4$.
Thus,
\begin{align*}
&    \|Q(y)-Q(.-\tfrac 12 \Delta)\|_{H^1(x<-T_c/2)}\\ &\leq  \|Q(y)\|_{H^1(x<-T_c/2)}+\|Q(.-\tfrac 12 \Delta)\|_{H^1(x<-T_c/2)}\leq K    e^{-\frac 12 c^{-\frac 1{100}}}.
\end{align*}

\medskip

\noindent 3. Position of the soliton $Q_c$ at $t=T_c$. We claim
    \begin{equation}\label{posQc1}
        \|Q_c(y_c)-b_{1,0}Q_c'(y_c)-Q_c(.+(1-c)T_c-b_{1,0})\|_{H^1}\leq K c^{\frac 1{p-1}+\frac 34}.
    \end{equation}
Indeed, for the $L^2$-norm, we have by a scaling argument
\begin{align*}
    \|Q_c-b_{1,0}Q_c'-Q_c(.-b_{1,0})\|_{L^2}& =c^{\frac 1{p-1}-\frac 14}
     \|Q-\sqrt{c}\,b_{1,0}Q'-Q(.-\sqrt{c}\,b_{1,0})\|_{L^2} \\
    &\leq K c^{ \frac 1{p-1}-\frac 14}(\sqrt{c}\, b_{1,0})^2 =K c^{ \frac 1{p-1}+\frac 34},
\end{align*}
and similarly for the estimate on the $x$ derivative.

Thus Proposition \ref{AVRILf} is proved.

\subsection{Extension of Proposition \ref{AVRILf} by scaling}
Let
$$ T_{c_1,c_2}=c_1^{-\frac 32}  T_c = c_1^{-\frac 32} \left(\frac {c_2}{c_1}\right)^{-\frac 12 -\frac 1{100}}.
$$
By a scaling argument, we have from Proposition \ref{AVRILf} the following

\begin{theorem}\label{corAVRIL}
Let $0<c_1<c_*(f)$ be such that \eqref{stable} holds. There exist $c_0(c_1)$ and $K_0(c_1)>0$, continuous in $c_1$
such that for any $0<c_2<c_0(c_1)$, there exist 
      function $v=v_{c_1,c_2}$ satisfying $v(0,x)=v(0,-x)$ and  
    such that  the following hold.
    \begin{enumerate}
        \item Approximate solution on $[-T_{c_1,c_2},T_{c_1,c_2}]$:  for $j= 0, 1,2$,
        \begin{equation}\label{AV0cor}
            \forall t\in [-T_{c_1,c_2},T_{c_1,c_2}],\quad   \|\partial_x^{j} S(t)\|_{L^2(\mathbb{R})}
                \leq K_0 c_2^{ \frac 2{p-1}+\frac 34 }.
        \end{equation}
        \item Closeness to the sum of two solitons for $t=\pm T_{c_1,c_2}$: there exist $\Delta_1$, $\Delta_2$ such that
\begin{equation}\label{3-16cor}\begin{split}
 &   \|v(T_{c_1,c_2})-Q_{c_1}(.-\tfrac 12\Delta_1)-Q_{c_2}(.+(c_1-c_2)T_{c_1,c_2}-\tfrac 12\Delta_2) \|_{H^1}
\leq K_0 c_2^{ \frac 2{p-1}+\frac 14 },\\&
    \|v(-T_{c_1,c_2})-Q_{c_1}(.+\tfrac 12\Delta_1)-Q_{c_2}(.-(c_1-c_2)T_{c_1,c_2}+\tfrac 12\Delta_2) \|_{H^1}
\leq K_0 c_2^{ \frac 2{p-1}+\frac 14 },
\end{split}\end{equation}
                where
        \begin{equation}\label{defDeltabiscor}
            \left|\Delta_1 - \left(\frac {c_2}{c_1}\right)^{\frac 1{p-1} -\frac 12} \delta_1 \right|\leq K c_2^{\frac 2{p-1}-\frac 12},\quad
        \delta_1=2\, \frac {\int Q_{c_1} \frac {d}{d c} {\int Q_{ c}}_{| c=c_1}}
        {\frac {d}{d c} \left({\int Q_{ c}^2}\right)_{| c=c_1}}.
        \end{equation}
        \item Closeness to the sum of two solitons: for all $t\in [-T_{c_1,c_2},T_{c_1,c_2}]$, there
        exists $y_1(t)$ such that
        \begin{equation}\label{OUBLI2}
        \| v(t,x)-Q_{c_1}(.-y_1(t))-Q_{c_2}(.-(c_2-c_1)t)\|_{H^1}\leq K_0 c^{\frac 1{p-1}}.
        \end{equation}
    \end{enumerate}
\end{theorem}

\noindent\emph{Proof of Theorem \ref{corAVRIL}.}
Fix a nonlinearity $f$ satisfying \eqref{surf}. Fix $0<c_1<c_*(f)$ such that \eqref{stable} holds.
Let
$$\widetilde f(\widetilde u)=\widetilde u^p + \widetilde f_1(\widetilde u)
\quad \text{where}\quad
\widetilde f_1(\widetilde u)=c_1^{-\frac p{p-1}} f_1(c_1^{\frac 1{p-1}} \widetilde u).$$
Then $u(t)$ is solution of \eqref{kdvf} if and only if
\begin{equation}\label{equiv}\widetilde u(t,x)= c_1^{-\frac 1{p-1}} u( c_1^{-\frac 32}t,c_1^{-\frac 12} x) \text{ is solution of }
  \partial_t \widetilde u + \partial_x (\partial_x^2 \widetilde u + \widetilde f(\widetilde u))=0.
\end{equation}

First, we observe that $\widetilde f$ satisfies assumption \eqref{surf}.
Second, for any $0<c<c_*(f)$, let $Q_c$ be the positive even solution of \eqref{ellipticf}.
For $0<\widetilde c=\frac c{c_1}<\frac {c_*(f)}{c_1}$, 
\begin{equation}\label{elliptictilde}
 \widetilde Q_{\widetilde c}(x)=c_1^{-\frac 1{p-1}}Q_c(c_1^{-\frac 12}x)  \quad \text{solves} \quad
\widetilde Q_{\widetilde c}''+\widetilde f(\widetilde Q_{\widetilde c})= \widetilde c\, \widetilde Q_{\widetilde c}.
\end{equation}
Thus, $c_*(\widetilde f)\geq \frac{c_*(f)}{c_1}>1$ (in fact, $c_*(\widetilde f)= \frac{c_*(f)}{c_1}$). Moreover,
for any $0<c<c_*(f)$, we have
\begin{equation}\label{scal2l1}\begin{split}
&
{\int Q_c^2}
= c_1^{\frac 2{p-1}-\frac 12} \int \widetilde Q_{\frac {c}{c_1}}^2,\quad
{\int Q_c}
= c_1^{\frac 1{p-1}-\frac 12} \int \widetilde Q_{\frac {c}{c_1}},\\
& \frac d{dc} {\int Q_c^2}_{| c=c_1}
= c_1^{\frac 2{p-1}-\frac 12}  \frac d{dc}\left( \int \widetilde Q_{\frac {c}{c_1}}^2\right)_{| c=c_1}
=c_1^{\frac 2{p-1}-\frac 32} \frac d{d\widetilde c}{\int \widetilde Q_{\widetilde c}^2}_{| \widetilde c=1},\\
&
\frac d{dc} {\int Q_c}_{| c=c_1}
=c_1^{\frac 1{p-1}-\frac 12} \frac d{dc}\left( \int \widetilde Q_{\frac {c}{c_1}}\right)_{| c=c_1}
=c_1^{\frac 1{p-1}-\frac 32} \frac d{d\widetilde c}{\int \widetilde Q_{\widetilde c}}_{| \widetilde c=1}.
\end{split}\end{equation}
In particular, $\frac d{dc} {\int Q_c^2}_{\big| c=c_1}>0$ is equivalent to 
$\frac d{d\widetilde c}{\int \widetilde Q_{\widetilde c}^2}_{\big| \widetilde c=1}>0$.

Let $c_0=\frac 14 c_0(\widetilde f)$, $K_0=K_0(\widetilde f)$, where $c_0(\widetilde f)$, $K_0(\widetilde f)$ are defined in Proposition \ref{AVRILf}
(these constants thus depend continuously upon $c_1$, see Remark after Proposition \ref{AVRILf}).
Let $0<c_2<c_0$, and let $c=\frac {c_2}{c_1}$. We consider $\widetilde v=\widetilde v_{1,c}$ as defined
in Proposition \ref{AVRILf} for the nonlinearity  $\widetilde f$ and
$\widetilde S=\partial_{t} \widetilde v + \partial_{x} (\partial_{x}^2 \widetilde v -\widetilde v + \widetilde f(\widetilde v))$.
From Proposition \ref{AVRILf}, we have
\begin{equation}\label{AV0pp}
            \forall t\in [-T_c,T_c],\quad   \|\partial_x^{j} \widetilde S(t)\|_{L^2(\mathbb{R})}
                \leq K_0 c^{ \frac 2{p-1}+\frac 34 }.
        \end{equation}
\begin{equation}\label{3-16pp} 
  \|\widetilde v(T_c)-\widetilde Q(.-\tfrac 12\widetilde \Delta)-\widetilde Q_c(.+(1-c)T_c-\tfrac 12\widetilde \Delta_c) \|_{H^1}\leq K_0 c^{ \frac 2{p-1}+\frac 14 },
\end{equation}
       \begin{equation}\label{defDeltabispp}
            \left|\widetilde \Delta - c^{\frac 1{p-1} -\frac 12} \widetilde \delta \right|\leq K c^{\frac 2{p-1}-\frac 12},
        \quad 
        \widetilde \delta=2 \, \frac {\int \widetilde Q \, \frac {d}{d\widetilde c} {\int \widetilde Q_{\widetilde c}}_{|\widetilde c=1}}
        {\frac {d}{d\widetilde c} \left({\int \widetilde Q_{\widetilde c}^2}\right)_{|\widetilde c=1}}.
        \end{equation}

Then, we set
\begin{equation}\label{defTccor}
v(t,x)=v_{c_1,c_2}(t,x)=c_1^{\frac 1{p-1}}  \widetilde v( c_1^{ \frac 32}t,c_1^{ \frac 12} x),
\end{equation}
\begin{equation}\label{defofScor}
        S(t,x)=\partial_{t} v + \partial_{x} (\partial_{x}^2 v - v + f(v)).
        \end{equation} 
Since $\partial_x^j S(t,x)=c_1^{\frac {3+j}2 +\frac 1{p-1}} \partial_x^j \widetilde S$, 
estimate \eqref{AV0pp} gives
$j=0,1,2 $, $\|\partial_x^{j} S(t)\|_{L^2(\mathbb{R})} \leq K  c_2^{ \frac 2{p-1}+\frac 34 }.$

From \eqref{3-16pp}
\begin{equation*} \|v(T_{c_1,c_2})-Q_{c_1}(.-\tfrac 12  c_1^{-\frac 12} \widetilde \Delta)-Q_{c_2}(.+(c_1-c_2) T_{c_1,c_2}-\tfrac 12  c_1^{-\frac 12}\widetilde \Delta_c
)\|_{H^1} \leq K  c_2^{\frac 2{p-1} +\frac 14}.
 \end{equation*}
Setting $\Delta_1=c_1^{-\frac 12} \widetilde \Delta $ and $\Delta_2= c_1^{-\frac 12} \widetilde \Delta_c$, by
\eqref{defDeltabispp} and \eqref{scal2l1}, we have
$$
\left|\Delta_1 - \left(\frac {c_2}{c_1}\right)^{\frac 1{p-1} -\frac 12} \delta_1 \right|\leq K c_2^{\frac 2{p-1}-\frac 12}
$$
$$
\delta_1 =   c_1^{ -\frac 12} \widetilde \delta
= 2\,  c_1^{ -\frac 12} \, \frac {\int \widetilde Q \frac {d}{d\widetilde c} {\int \widetilde Q_{\widetilde c}}_{|\widetilde c=1}}
        {\frac {d}{d\widetilde c} \left({\int \widetilde Q_{\widetilde c}^2}\right)_{|\widetilde c=1}}
= 2\, \frac {\int Q_{c_1} \frac {d}{d c} {\int Q_{c}}_{| c=c_1}}
        {\frac {d}{d c} \left({\int Q_{c}^2}\right)_{| c=c_1}}.
$$
Estimate \eqref{OUBLI2} follows from \eqref{OUBLI1}.

\section{Preliminary results for stability of the $2$-soliton structure}
This section is similar to Section 4 in \cite{MMcol1}. 

\subsection{Dynamic stability in the interaction region}
\begin{proposition}[Exact solution close to the approximate solution $v$]\label{INTERACT4col}
    Let $0<c_1<c_*(f)$ be such that \eqref{stable} holds. There exist $c_0(c_1)$ and $K_0(c_1)>0$, continuous in $c_1$
such that for any $0<c_2<c_0(c_1)$, the following holds. Let $v=v_{c_1,c_2}$ be  defined in Theorem
\ref{corAVRIL}.
    Suppose that for some $\theta>\frac 1{p-1}$,
    for some $T_0\in [-T_{c_1,c_2},T_{c_1,c_2}]$,
    \begin{equation}\label{hypINT}
        \| u(T_0) - v(T_0) \|_{H^1(\mathbb{R})}\leq  c_2^{\theta},
    \end{equation}
where  $u(t)$ is an $H^1$ solution of \eqref{kdvf}. 
    Then, $u(t)$ is global and 
    there exists   $\rho(t)$ such that, for all $t\in [-T_{c_1,c_2},T_{c_1,c_2}],$
    \begin{equation}\label{INT41}
        \|u(t)-v(t,.-\rho(t)) \|_{H^1}+ |\rho'(t)-c_1|\leq 
        K_0 \left(c_2^{\theta}+ c_2^{\frac 2{p-1} +\frac 14 -\frac 1{100}}\right).
    \end{equation}
\end{proposition}

The fact that $u(t)$ is global follows from the stability of $Q_{c_1}$.

\medskip

\noindent\emph{Sketch  of the proof of Proposition \ref{INTERACT4col}.}
  The   proof is similar to the one of Proposition 4.1 in \cite{MMcol1}.
For the sake of simplicity, we give a sketch of the proof in the special case $c_1=1$ and $c_2=c$ small, i.e.
we work in the context of Proposition \ref{AVRILf}. The general case follows by the same scaling argument as in Section 2.5. In view of \eqref{INT41}, we may assume that 
\begin{equation}\label{theta}
\theta\leq \frac 2{p-1} +\frac 14 -\frac 1{100}.
\end{equation}
We  prove the result on $[T_0,T_c]$. By using the transformation $x\to -x$, $t\to -t$, the proof is the same on $[-T_c,T_0]$.

Let $K^*>1$ be a constant to be fixed later. Since $\|u(T_0)-v(T_0)\|_{H^1}\leq c^\theta$, by continuity in time in $H^1(\mathbb{R})$, there exists $T_0<T^*\leq T_c$ such that
\begin{equation*}
    T^*=\sup\left\{T\in [T_0,T_c] \text{ s.t. $\forall t\in [T_0,T]$, $\exists r(t)\in \mathbb{R}$ with }
    \|u(t){-}v(t,.{-}r(t))\|_{H^1}\leq K^* c^{\theta} \right\}.
\end{equation*}
The objective is to prove that $T^*=T_c$ for $K^*$ large. For this, we argue by contradiction, assuming that $T^*<T_c$ and reaching a    contradiction with the definition of $T^*$ by proving independent estimates on $\|u(t)-v(t,.-r)\|_{H^1}$ on $[T_0,T^*]$.

We claim (see Lemma 4.1 in  \cite{MMcol1}).
\begin{claim}\label{DEFZ} Assume that $0<c<c(K^*)$ small enough.
There exists a unique $C^1$ function $\rho(t)$ such that, for all $t\in [T_0,T^*]$,
\begin{equation}\label{defz}
z(t,x)=u(t,x+\rho(t))-v(t,x) \quad \text{satisfies}\quad 
 \int z(t,x) Q'(y) dx=0.
\end{equation}
Moreover, we have,  for all $t\in [T_0,T^*]$,
\begin{equation}\label{TRANS3}
              |\rho(T_0)|+\|z(T_0)\|_{H^1}\leq K c^{\theta}, \ \|z(t)\|_{H^1}\leq  2 K^* c^{\theta},
\end{equation}
\begin{equation}\label{eqz}
\partial_t z +\partial_x (\partial_x^2 z -z + f(z+v)-f(v))= -S(t) + (\rho'(t)-c_1) \partial_x (v+z).
\end{equation}
\begin{equation}\label{TRANS3bis}
 |\rho'(t)-1|\leq K \|z(t)\|_{H^1}+K \|S(t)\|_{H^1},
\end{equation}
\end{claim}
Recall that the existence, uniqueness and regularity of $\rho(t)$ is obtained by
a standard use of the Implicit Function Theorem applied to $u(t)$ at each fixed time $t$.
Estimate \eqref{TRANS3bis} is obtained by equation \eqref{eqz}.

\medskip

\emph{Step 1.} Energy estimates on $z(t)$.
We extend to the case of the general power nonlineartity the
definition given in \cite{MMcol1}  of the energy functional for $z(t)$:
$$\mathcal{F}(t)=\frac 12 \int \left((\partial_x z)^2 + (1+\alpha'(y_c))z^2\right)
- \int (F(v+z) - F(v) - f(v) z).$$

\begin{lemma}[Coercivity of $\mathcal{F}$]\label{varF}
Assume that $0<c<c(K^*)$ small enough.
    There exists $K>0$ (independent of $K^*$ and $c$) such that
    \begin{description}
    \item{\rm (i)} Coercivity of $\mathcal{F}$ under orthogonality conditions:
    \begin{equation}\label{posf1}
       \forall t\in [T_0,T^*],\quad\|z(t)\|_{H^1}^2\leq K \mathcal{F}(t) + K \left|\int z(t) Q(y)\right|^2.
    \end{equation}
    \item{\rm (ii)} Control of the direction $Q$:
    \begin{equation}\label{Qdir1} \forall t\in [T_0,T^*],\quad
        \left|\int  z(t)Q(y)\right|\leq  K c^\theta + K c^{\frac 1{p-1} -\frac 14} 
        \|z(t)\|_{L^2}+ K\|z(t)\|_{L^2}^2.
\end{equation}
\item{\rm (iii)} Control of the variation of the energy fonctional:
\begin{equation}\label{varF1}
        \mathcal{F}(T^*)-\mathcal{F}(T_0)\leq  K c^{2 \theta}\left((K^*)^2(1+K^*)c^{\frac 1{2(p-1)} -\frac 18} + K^*\right).
    \end{equation}
\end{description}\end{lemma}

\noindent\emph{Proof of Lemma \ref{varF}.}
(i) For this property, see proof of  Claim 4.2 in Appendix D of \cite{MMcol1}.
Recall that the proof of such property is related to assumption \eqref{stable} (nonlinear stability of 
$Q$) and to the choice
of $\rho(t)$ in Claim \ref{DEFZ}.

\medskip

(ii) This estimate follows from the conservation of $\int u^2(t)$ and a similar
approximate conservation for $v(t)$. Indeed, we have
$
  |\frac 12 \frac d {dt} \int v^2 |= | \int S(t,x) v(t,x) dx |
    \leq K \|S(t)\|_{L^2}
$
from the equation of $v(t)$ (see \cite{MMcol1} for more details).

\medskip

(iii) The computations of the proof of Lemma 4.3 in \cite{MMcol1} are extended as follows:
$$\mathcal{F}'(t)= \mathbf{F}_1+\mathbf{F}_2+\mathbf{F}_3,$$
where
$$\mathbf{F}_1=\int \partial_t z (- \partial_x^2 z + z - (f(v+z)-f(v))),
\quad \mathbf{F}_2= \int \partial_t z \, \alpha'(y_c) z,$$
$$\mathbf{F}_3= \int \left\{ \frac 12 (1-c) \alpha''(y_c) z^2 
-\partial_t v \left(f(v+z)-f(v)-z f'(v)\right).
\right\}$$
Then, we have, for $m_0=\min\left(\frac 2{p-1},\frac 1{p-1} +\frac 12\right)$,
\begin{equation}\label{assum1}
\left|\mathbf{F}_1 + (\rho'(t)-1) \int \alpha'(y_c) Q'(y) z \right|
\leq K c^{\frac 1{p-1}+\frac 14} \|z(t)\|_{L^2}^2 + K \|z(t)\|_{L^2} 
(\|\partial_x^2 S(t)\|_{L^2}+\|S(t)\|_{L^2}),
\end{equation}
\begin{equation}\label{assum2}\begin{split}
&\left| \mathbf{F}_2 - (\rho'(t)-1) \int \alpha'(y_c) Q'(y) z + \frac 12
\int \alpha'(y_c) Q'(y) f''(Q(y)) z^2  \right|\\ &\leq 
K  \|z(t)\|_{H^1}^2 \left(c^{m_0}+ c^{\frac 1{p-1}} \|z(t)\|_{H^1}\right)
+ K \|z(t)\|_{H^1} (\|\partial_x^2 S(t)\|_{L^2}+\|S(t)\|_{L^2}),
\end{split}\end{equation}
\begin{equation}\label{assum3}
\left| \mathbf{F}_3 - \frac 12
\int \alpha'(y_c)  Q'(y) f''(Q(y))  z^2 \right|\leq
K c^{m_0} \|z(t)\|_{H^1}^2 + K c^{\frac 1{p-1}} \|z(t)\|_{H^1}^3. 
\end{equation}
Estimates \eqref{assum1}--\eqref{assum3} are obtained exactly as in 
\cite{MMcol1}. Now, we conclude the proof of Lemma \ref{varF}.

From the cancellations of the main terms of $\mathbf{F}_1$, $\mathbf{F}_2$ and $\mathbf{F}_3$, and then from  \eqref{TRANS3} and Theorem \eqref{corAVRIL}, \eqref{AV0cor},   we  get
\begin{equation*}\begin{split}
|\mathcal{F}'(t)| &
\leq           K \|z(t)\|_{H^1}^2 \left(c^{\frac 1{p-1}+\frac 14}+c^{\frac 1{p-1}} \|z(t)\|_{H^1}\right)+K \|z(t)\|_{H^1}\left(\|\partial_x^2 S(t)\|_{L^2}+\|S(t)\|_{L^2}\right)\\
& \leq  K c^{2 \theta} \left[ (K^*)^2  ( c^{\frac 1{p-1}+\frac 14}+ K^*c^{\frac 1{p-1}+\theta})+ K^* c^{ \frac 2{p-1}+\frac 34- \theta }\right].
\end{split}
\end{equation*}
Integrating on the time interval $[T_0,T^*]$, since $T^*-T_0 \leq  2T_c= 2 c^{\frac 12 +\frac 1{100}}$, 
and $\theta>\frac 1{p-1}>\frac 14$, we obtain
\begin{equation*}
    |\mathcal{F}(T^*)-\mathcal{F}(T_0)|\leq   K c^{2 \theta}\left((K^*)^2(1+K^*)   c^{\frac 1{p-1}-\frac 14-\frac 1{100}} + K^* c^{ \frac 2{p-1}+\frac 14-\frac 1{100}- \theta }\right].    \end{equation*}
Note that by \eqref{theta}, we have $ \frac 2{p-1}+\frac 14-\frac 1{100}- \theta\geq 0$ and 
$\frac 1{p-1}-\frac 14-\frac 1{100}  \geq \frac {1}{2(p-1)} -\frac 18 >0$,
since $\frac {1}{2(p-1)} \geq \frac 16 \geq \frac 18+\frac 1{100}$.
Thus, Lemma \ref{varF} is proved.

\medskip

\emph{Step 2.} Conclusion of the proof.
By  \eqref{Qdir1}, we have
\begin{equation*}
    \left|\int  z(T^*) Q(y) \right|\leq  K c^\theta + K c^{\frac 1{p-1}-\frac 14 } \|z(T^*)\|_{L^2}
    + \|z(T^*)\|_{L^2}^2,
\end{equation*}
and thus by \eqref{posf1},
\begin{equation*}
    \|z(T^*)\|_{H^1}^2\leq K \mathcal{F}(T^*) + K (c^\theta +   c^{\frac 1{p-1}-\frac 14} \|z(T^*)\|_{L^2}+ \|z(
    T^*)\|_{L^2}^2)^2.
\end{equation*}
Since $\frac 1{p-1}-\frac 14>0$, it follows that for $c$ small enough,
\begin{equation*}
    \|z(T^*)\|_{H^1}^2\leq (K+1) \mathcal{F}(T^*) + K c^{2 \theta}.
\end{equation*}

Next, by \eqref{varF1} and $|\mathcal{F}(T_0)|\leq K c^{2 \theta}$, we obtain
\begin{equation*}
    \|z(T^*)\|_{H^1}^2 
\leq (K+1) (\mathcal{F}(T^*)-\mathcal{F}(T_0))+K c^{2\theta}
\leq K_1 c^{2 \theta}\left((K^*)^2(1+K^*) c^{\frac 1{2(p-1)} -\frac 18} + K^*+1\right),
\end{equation*}
where $K_1$ is independent of $c$ and $K^*$.
Choose $c_*=c_*(K^*)$ such that  
\begin{equation*}
    (K^*)^2(1+K^*) c_*^{\frac 1{2(p-1)} -\frac 18}<1.
\end{equation*}
Then, for $0<c<c_*$,
\begin{equation*}
    \|z(T^*)\|_{H^1}^2 \leq K_1 c^{2 \theta}\left(2+ K^*\right).
\end{equation*}
Next, fix $K^*$ such that $K_1 (2+K^*)<\frac 12 (K^*)^2$. Then 
\begin{equation*}
    \|z(T^*)\|_{H^1}^2 \leq \frac 12 (K^*)^2 c^{2 \theta}.
\end{equation*}
This contradict the definition of $T^*$, thus proving that $T^*=T_c$.
Thus estimate \eqref{INT41} is proved on $[T_0,T_c]$.

\subsection{Stability and asymptotic stability for large time}
 
In this section, we consider the stability of the $2$-soliton structure after the collision. This question
has been considered in \cite{MMas1}, \cite{MMas2}. See also \cite{MM1}, \cite{MMT}, \cite{yvanSIAM}.
We recall the following.

\begin{proposition}[Stability and asymptotic stability \cite{MMas1}, \cite{MMas2}]\label{ASYMPTOTIC}
Let $0<c_1<c_*(f)$ be such that \eqref{stable} holds. There exist $c_0(c_1)$ and $K_0(c_1)>0$, continuous in $c_1$
such that for any $0<c_2<c_0(c_1)$ and for any $\omega>0$, the following hold.
Let $u(t)$ be an $H^1$ solution of \eqref{kdvf} such that
for some $t_1\in \mathbb{R}$ and $\frac 12 {T_{c_1,c_2}}\leq X_0\leq \frac 32{T_{c_1,c_2}}$,
\begin{equation}
\label{D25}
\|u(t_1)-Q_{c_1}-Q_{c_2}(.+X_0)\|_{H^1}\le c_2^{\omega+\frac 1{p-1}+\frac 14}.
\end{equation}
Then, there exist    $C^1$  functions $\rho_1(t)$, $\rho_2(t)$ defined on $[t_1,+\infty)$  such that
\begin{enumerate}
\item Stability.
\begin{equation}\label{huit}
\sup_{t\ge t_1} \| u(t)-Q_{c_1}(.-\rho_1(t))-Q_{c_2}(.-\rho_2(t))  \|_{H^1 }
\le K   c^{\omega+\frac 1{p-1}-\frac 14},
\end{equation}
\begin{equation}\label{suppl}\begin{split}
& \forall t\ge t_1,\ \tfrac 12 c_1\leq (\rho_1-\rho_2)'(t) \leq \tfrac 32 c_1,
\\ &|\rho_1(t_1)|\leq K c_2^{\omega+\frac 1{p-1}+\frac 14},\quad
|\rho_2(t_1)-X_0|\leq K c_2^{\omega}.
\end{split}\end{equation}
\item Convergence of $u(t)$.
There exist $c_1^+, c_2^+ >0$  such that
\begin{equation}\label{neuf}
\lim_{t\rightarrow +\infty}\|u(t)-Q_{c_1^+}(x-\rho_1(t))-Q_{c_2^+}(x-\rho_2(t))\|
_{H^1(x> \frac{c_2t}{10})}=0.
\end{equation}
\begin{equation}\label{sept}
\left|\frac{c_1^+}{c_1}-1\right|\leq   K  c^{\omega+\frac 1{p-1}+\frac 14 },\quad
 \left|\frac {c_2^+}{c_2} - 1\right|\le K  c^{\omega}.
\end{equation}
\end{enumerate}
\end{proposition}
The proof of 
Proposition \ref{ASYMPTOTIC} is based on energy arguments, monotonicity results on local energy quantities,
and a Virial argument on the linearized problem around solitons.

The loss of $\frac 12$ in the exponent between \eqref{D25} and \eqref{huit}
is due to the fact that the natural norm to study the stability of $Q_{c_2}$ is
not $\|.\|_{H^1}$ but $\|\partial_x(.)\|_{L^2}+c^{\frac 12}\|.\|_{L^2}$.

\subsection{Monotonicity results}
Recall a more precise decomposition of $u(t)$ used in the proof of Proposition \ref{ASYMPTOTIC} in
\cite{MMas1}, \cite{MMas2}.

\begin{claim}[Decomposition of the solution]\label{LEMMEB1}
Under the assumptions of Proposition \ref{ASYMPTOTIC}, 
 there exist $C^1$  functions $\rho_1(t)$, $\rho_2(t)$, $c_1(t),$ $c_2(t)$,
defined on $[t_1,+\infty)$,  such that the function $\eta(t)$ defined by
\begin{equation*}
\eta(t,x)=u(t,x)-R_1(t,x)-R_2(t,x),
\end{equation*}
where for $ j=1,2,$
$
R_j(t,x)=Q_{c_j(t)}(x-\rho_j(t)),$
satisfies for all $t\geq t_1,$
\begin{eqnarray}&&
\int R_j(t)\eta(t) =\int (x-\rho_j(t))  R_j(t)\eta(t)=0,\quad j=1,2,\label{dix}\\ &&
\|\eta(t)\|_{H^1}+\left|\frac{c_1(t)}{c_1}-1\right|
+c_2^{\frac 1{p-1}-\frac 14}\left|\frac {c_2(t)}{c_2} -1\right|
\le K   c_2^{\omega+\frac 1{p-1}-\frac 14},\label{onze}
\end{eqnarray}
\end{claim}

Now, we recall some   monotonicity results for two localized quantities defined in $\eta(t)$.
Define 
\begin{equation}\label{surphiff} 
 \psi(x)=\tfrac {2}{\pi}\arctan(\exp(- \tfrac 14 x )),
\end{equation}
\begin{equation}\label{defgj}
g_j(t)= \int(\eta_x^2 + c_j \eta^2) (t,x) e^{-\frac 14 \sqrt{c_j} |x-\rho_j(t)|} dx,\quad j=1,2.
\end{equation}
For  $0\leq t_0\leq t$, $x_0\geq 0$, $j=1,2$, let
\begin{equation*}
\begin{split}
& \mathcal{M}_j(t)=\int \eta^2 \psi_j ,\\
&\mathcal{E}_j(t)=\int  \left[\frac 12 \eta_x^2
-\left( F(R_1{+}R_2{+}\eta)  {-}(f(R_1) +f(R_2)) \eta {-}F(R_1{+}R_2)\right)\right] \psi_j,
\end{split}
\end{equation*}
\begin{equation*}\begin{split}
\hbox{where}\quad & \psi_1(x)=\psi(\sqrt{c_1} \widetilde x_1),
\quad \widetilde x_1= x-\rho_1(t)+x_0+\tfrac {c_1}2 (t-t_0), \\
& \psi_2(x)=\psi(\sqrt{c_2}\widetilde x_2),\quad
\widetilde x_2= x-\rho_2(t)+x_0+ \tfrac {c_2}2 (t-t_0).
\end{split}\end{equation*}
 \begin{claim}[Monotonicity results in $\eta(t)$]\label{NEW}
Let $x_0>0$, $t_0>0$.
For all $t\geq t_0 $,
\begin{equation*}
\begin{split}
	& \frac d{dt}\left( \int Q^2_{c_1(t)} + \mathcal{M}_1(t) \right) 
	   \leq K e^{-\frac {\sqrt{c_1}} {16} (c_1(t-t_0)+x_0)} g_1(t) + Ke^{-\frac 1 {32} c_1 \sqrt{c_2} (t+T_{c_1,c_2})},\\
	& \frac d{dt}\left(2E(Q_{c_1(t)})+ 2 \mathcal{E}_1(t) 
	+\frac {c_1}{100} \left(\int Q^2_{c_1(t)}+ \mathcal{M}_1(t) \right)
	\right)
	   \\ &\qquad \leq K e^{-\frac 1{16}\sqrt{c_1}(c_1(t-t_0)+x_0)} g_1(t)+ 
	   Ke^{-\frac 1 {32} c_1 \sqrt{c_2} (t+T_{c_1,c_2})}.\\
	   & \frac d{dt}\left( \int Q^2_{c_1(t)} + \int Q^2_{c_2(t)}  + \mathcal{M}_2(t) \right) 
	   \leq K e^{-\frac {c_2\sqrt{c_2}}{16}    (t-t_0)} e^{-\frac {\sqrt{c_2} }{16}  x_0} \sqrt{c_2}\, g_2(t) + Ke^{-\frac 1 {32} c_1 \sqrt{c_2 } (t+T_{c_1,c_2})},\\
	& \frac d{dt}\left(2E(Q_{c_1(t)}) {+}2E(Q_{c_2(t)}) + 2 \mathcal{E}_2(t) 
	+\frac {c_2}{100} \left( \int Q^2_{c_1(t)}{+} \int Q^2_{c_2(t)} \right)+ \mathcal{M}_2(t) \right)
	 \\ &\qquad
	   \leq K e^{-\frac {c_2\sqrt{c_2}}{16}   (t-t_0)  } e^{-\frac {c_2 }{16}  x_0} c_2^{\frac 32} g_2(t)
	   + Ke^{-\frac 1 {32} c_1 \sqrt{c_2} (t+T_{c_1,c_2})}.
\end{split}
\end{equation*}
\end{claim}
Claim \ref{NEW} is proved in \cite{MMas2} for the power case. The proof is exactly
the same for a nonlinearity $f(u)$ satisfying \eqref{surf}.

\section{Proof of the main Theorems}

\subsection{Proof of Theorem \ref{PURE}}
Let $0<c_1<c_*(f)$ such that \eqref{stable} holds and $c_2>0$ small enough.
Let $u(t)$ be the unique solution of \eqref{kdvf} such that
(see Theorem 1 and Remark 2 in \cite{Martel})
$$
\lim_{t\to -\infty} \|u(t)-Q_{c_1}(x-c_1 t)-Q_{c_2}(x-c_2 t)\|_{H^1}=0.
$$

\medskip

\noindent\emph{1. Behavior at $-T_{c_1,c_2}$.}
We claim that
\begin{equation}\label{ETI1}
\forall t<-\frac 1{32} T_{c_1,c_2},\quad 
\|u(t)-Q_{c_1}(.-c_1 t)-Q_{c_2}(.-c_2 t)\|_{H^1}
\leq K e^{\frac 14 \sqrt{c_2}(c_1-c_2)t}.
\end{equation}
This is a consequence of the proof of existence of $u(t)$ in \cite{Martel}. 
See Proposition 5.1 in \cite{MMcol1} for a proof in the power case.

Now, let $\Delta_1$, $\Delta_2$ be defined in Theorem \ref{corAVRIL} and
$$
T_{c_1,c_2}^- = T_{c_1,c_2} + \frac 12 \frac {\Delta_1-\Delta_2}{c_1-c_2},\quad
a=\frac 12 \Delta_1 - T_{c_1,c_2}^-.
$$
Since $|\Delta_1|\leq K c^{-\frac 16}$ and $\Delta_2$ is independent of $c$,
we have
$
-T_{c_1,c_2}^-\leq -\frac 1{32} T_{c_1,c_2}$, and thus, for $c_2$ small enough:
\begin{equation*}
\|u(-T_{c_1,c_2}^-,.+a)-Q_{c_1}(.+\tfrac {\Delta_1}2)-Q_{c_2}(.-(c_1-c_2)T_{c_1,c_2}+\tfrac {\Delta_2}2)\|_{H^1} 
\leq K e^{- \frac 14 \sqrt{c_2}(c_1-c_2)T_{c_1,c_2}^-}\leq K c_2^{10}.
\end{equation*}
Let $\widetilde u(t,x)=u(t-T_{c_1,c_2}+T^-_{c_1,c_2},x-a)$. Then $\widetilde u(t,x)$ is also solution of 
\eqref{kdvf} and satisfies
\begin{equation}\label{atmoinstc}
\|\widetilde u(-T_{c_1,c_2})-Q_{c_1}(.+\tfrac {\Delta_1}2)
-Q_{c_2}(.-(c_1-c_2)T_{c_1,c_2}+\tfrac {\Delta_2}2)\|_{H^1} \leq  K c_2^{10}.
\end{equation}
In what follows, we work with $\widetilde u(t)$ satisfying \eqref{atmoinstc} and we denote $\widetilde u$ by $u$.

\medskip

\noindent\emph{2. Behavior at $+T_{c_1,c_2}$.}  Now, consider $v=v_{c_1,c_2}$ constructed in Theorem \ref{corAVRIL} (possibly taking a smaller $c_2$).
By \eqref{3-16cor} and \eqref{atmoinstc}, we have
$$
\|u(-T_{c_1,c_2})-v(-T_{c_1,c_2})\|_{H^1} \leq K c_2^{\frac 2{p-1}+\frac 14}.
$$
Applying 
Proposition \ref{INTERACT4col} with 
$$T_0=-T_{c_1,c_2},\quad \theta =\frac 2{p-1} +\frac 14,$$
it follows that there exists
a function $\rho(t)$ such that
$$
\forall t\in [-T_{c_1,c_2},T_{c_1,c_2}],\quad
\|u(t)-v(t,.-\rho(t))\|_{H^1}\leq K c^{\frac 2{p-1}+\frac 14-\frac 1{100}}.
$$
In particular, by \eqref{3-16cor}, 
for some $a_-$, $b_-$ such that $\frac 12 T_{c_1,c_2} <a_--b_-<2 T_{c_1,c_2}$,
\begin{equation}\label{atplustc}
\|u(T_{c_1,c_2})-Q_{c_1}(.-a_-)-Q_{c_2}(.-b_-)\|_{H^1}
\leq K c^{\frac 2{p-1}+\frac 14-\frac 1{100}}.
\end{equation}

\smallskip

\noindent\emph{3. Behavior as $t\to +\infty$.}
From \eqref{atplustc}, it follows that we can apply Proposition \ref{ASYMPTOTIC} to
$u(t)$ for $t\geq T_{c_1,c_2}$, with $$\omega=\frac 1{p-1} -\frac 1{100}.$$
It follows that there exist $\rho_1(t)$, $\rho_2(t)$, $c_1^+$, $c_2^+$ so that
\begin{align}
& w^+(t,x)=u(t,x)
-Q_{c_1^+}(x-\rho_1(t))-Q_{c_2^+}(x-\rho_2(t))\quad \text{satisfies}\\
& \sup_{t\geq T_{c_1,c_2}} \|w^+(t)\|_{H^1} \leq  K c_2^{\frac 2{p-1}-\frac 14-\frac 1{100}},
\quad \lim_{t\to +\infty}
\|w^+(t)\|_{H^1(x>\frac {c_2}{10} t)} =0,\label{ETI5}\\
&
|c_1^+-c_1|\leq K c_2^{\frac 2{p-1}+\frac 14-\frac 1{100}},\quad
|c_2^+-c_2|\leq K c_2^{1+\frac 1{p-1}-\frac 1{100}}.\label{ETI3}
\end{align}

\smallskip

\noindent\emph{4. Estimates on $c_1^+-c_1$ and $c_2^+-c_2$.}
By \eqref{ETI1} and conservation of the $L^2$ norm, we have
$$
M_0=\int u^2(t)=\int Q_{c_1}^2+ \int Q_{c_2}^2.
$$
By the definition of $w^+(t)$, we have
$$
\forall t,\quad
M_0= \int Q_{c_1^+}^2 + \int Q^2_{c_2^+} + \int (w^+)^2(t)+
2 \int  w^+(t)(Q_{c_1^+}+Q_{c_2^+}) + 2\int Q_{c_1^+} Q_{c_2^+}.
$$
Thus, by \eqref{ETI5}, passing to the limit as $t\to +\infty$, we obtain
$M^+=\lim_{t\to +\infty} \int (w^+)^2(t)$ exists and
\begin{equation}\label{unun}
M^+
=\int Q_{c_1}^2 +\int Q_{c_2}^2 -\int Q_{c_1^+}^2 -\int Q_{c_2^+}^2, Ê
\end{equation}
Similarly, using the conservation of energy, 
$E^+=\lim_{t\to +\infty} E(w^+(t))$ exists and
\begin{equation}\label{deuxdeux}
 E^+ =E(Q_{c_1})+E(Q_{c_2})-E(Q_{c_1^+})-E(Q_{c_2^+}).
\end{equation}
By \eqref{ETI5}, we have $\|w^+(t)\|_{L^\infty}^{p-1}\leq K \|w^+(t)\|^{p-1}_{H^1}
\leq K c_2^{\frac 98}$, for $t$ large enough. Thus,
\begin{align*}
E(w^+(t)) & = \frac 12 \int (w_x^+)^2(t) - \int F(w^+(t))
\geq \frac 12 \int (w_x^+)^2(t) - K \|w^+(t)\|_{L^\infty}^{p-1} \int (w^+)^2(t)\\
& \geq \frac 12 \int (w_x^+)^2(t) - K \|w^+(t)\|_{L^\infty}^{p-1} \int (w^+)^2(t) \geq \frac 12 \int (w_x^+)^2(t) - K c_2^{\frac 98} \int (w^+)^2(t).
\end{align*}
Passing to the limit $t\to +\infty$, we obtain \eqref{th02bis}.

If $\limsup_{t\to +\infty} \|w^+_x(t)\|_{L^2} + \|w^+(t)\|_{L^2} = 0$,
then $w^+(t)\to 0$ in $H^1$ as $t\to +\infty$, and $u(t)$ is a pure two soliton solution 
at $+\infty$, $c_1^+=c_1$ and $c_2^+=c_2$ so that \eqref{th03}--\eqref{th04} hold.  

Assume now that $\limsup_{t\to +\infty} \|w^+_x(t)\|_{L^2} + \|w^+(t)\|_{L^2} >0$, so that
$E^++\tfrac 12 c_2 M^+>0$.
Recall that (\cite{We1}) by assumption \eqref{stable},
\begin{equation}\label{quatreet}
\frac d{dc} E(Q_c) = - \frac 12 c \frac d{dc} \int Q_c^2<0,
\quad \text{for $c=c_1$ and  $c=c_2$.}
\end{equation}
Let
$\bar c_2$ be such that $\bar c_2\left( \int Q^2_{c_2}-\int Q^2_{c_2^+}\right) = 
2 (E(Q_{c_2})-E(Q_{c_2}))$. Then, by \eqref{ETI3} and \eqref{quatreet} on $c_2$
we have $|\frac {\bar c_2}{c_2}-1|\leq \frac 14$.
Multiplying \eqref{unun} by $\bar c_2$ and summing \eqref{deuxdeux}, we find:
$$
E^++\tfrac {\bar c_2}2 M^+
= E(Q_{c_1})-E(Q_{c_1^+}) + \tfrac {\bar c_2} 2 \Big(\int Q_{c_1}^2-\int Q_{c_1^+}^2\Big).
$$
Using \eqref{ETI3} and \eqref{quatreet} on $c_1$, we find
\begin{equation}
            \frac 1{K} (2 E^+ + c_2 M^+) \leq \frac {c_1^+}{c_1} -1
            \leq {K}   (2 E^+ + c_2 M^+),
\end{equation}
Let $\bar c_1$ be such that $\bar c_1\left( \int Q_{c_1}^2-\int Q_{c_1^+}^2\right) = 
2 (E(Q_{c_1})-E(Q_{c_1}))$.
Arguing similarly, we have $|\bar c_1-c_1|\leq \frac 14 c_1$ and
$$
E^++\tfrac {\bar c_1}2 M^+
= E(Q_{c_2})-E(Q_{c_2^+}) + \tfrac {\bar c_1} 2 \Big(\int Q_{c_2}^2-\int Q_{c_2^+}^2\Big).
$$
 By \eqref{surf}, since $c_2$ is small, we have
 $\frac d{dc} {\int Q^2_c}_{|c=c_2}\sim (\frac 2{p-1} -\frac 12) c_2^{\frac 2{p-1} -\frac 32},$
and thus
\begin{equation}
        \frac 1K c_2^{\frac 2{p-1}-\frac 12} (2 E^+ + c_1 M^+) 
          \leq 1-\frac {c_2^+}{c_2} 
            \leq {K} c_2^{\frac 2{p-1}-\frac 12} (2 E^+ + c_1 M^+) .
\end{equation}
This concludes the proof of Theorem \ref{PURE}.

\subsection{Proof of existence. Theorem \ref{EXISTf}}
For $0<c_1<c_*(f)$ such that \eqref{stable} holds and $c_2>0$ small enough, we
denote by  $u_{c_1,c_2}(t)$ the global solution of 
\begin{equation}\label{defv12}
\partial_t u +\partial_x(\partial_x^2 u +f(u))=0,\quad
u(0,x)=v_{c_1,c_2}(0,x),
\end{equation}
where $v_{c_1,c_2}(t)$ is the approximate solution 
constructed in Theorem \ref{corAVRIL} (note that $u_{c_1,c_2}(t)$ is global by stability of $Q_{c_1}$).
By the parity property of $x\mapsto v_{c_1,c_2}(0,x)$ and since   equation \eqref{kdvf} is invariant
under the transformation $x\to -x$, $t\to -t$, the solution $u_{c_1,c_2}(t)$ has the following
symmetry: 
\begin{equation}\label{SYMM}
u_{c_1,c_2}(t,x)=u_{c_1,c_2} (-t,-x).
\end{equation}
Thus, we shall only study $u_{c_1,c_2}(t)$ for $t\geq 0$.
We claim the following concerning $u_{c_1,c_2}(t)$.

\begin{proposition}
\label{DANSTHM12}
Let $0<c_1<c_*(f)$ be such that \eqref{stable} holds.
There exist $c_0(c_1)>0$ and $K_0(c_1)>0$, continuous in $c_1$ such that 
for any $0<c_2<c_0(c_1)$,
there exist $0<c_2^+(c_1,c_2)<c_1^+(c_1,c_2)<c_*(f)$,
 and $\rho_1(t;c_1,c_2)$, $\rho_2^+(t;c_1,c_2)\in \mathbb{R}$, such that the following hold for
$$
w^+_{c_1,c_2}(t,x)=u_{c_1,c_2} (t,x)-Q_{c_1^+ }(x-\rho_1(t) )-Q_{c_2^+ }(x-\rho_2(t) ).
$$
\begin{enumerate}
\item Asymptotic behavior:
\begin{equation}\label{thm12I}
\lim_{t\rightarrow +\infty}\| w^+_{c_1,c_2}(t)\|
_{H^1(x> c_2 t/10)}=0.
\end{equation}
\begin{equation}\label{thm12IV} 
\text{for $t$ large},\quad 
 \|w^+_{c_1,c_2}(t)\|_{H^1}\leq K_0    c_2^{\frac{2}{p-1}+\frac 14 -\frac 1{100}},
\end{equation}
\begin{equation}\label{thm12II}
\left|\frac{c_1^+}{c_1}  - 1\right|\leq K_0 c_2^{\frac 2{p-1}+\frac 14 -\frac 1{100}}  ,\quad
 \left|\frac {c_2^+ }{c_2} - 1\right|\le K_0 c_2^{\frac 1{p-1}-\frac 1{100}},\end{equation}
 \begin{equation}\label{thm12III}
|\rho_1(T_{c_1,c_2} ){-}(c_1 T_{c_1,c_2} {+}\tfrac 12 \Delta_1) | \leq  K c_2^{\frac 2{p-1}-\frac 12},\quad 
 | \rho_2(T ){-}(c_2 T_{c_1,c_2} {+}\tfrac 12 {\Delta_2})  |\leq K c_2^{\frac 1{p-1}-\frac 1{100}},
\end{equation}
where $\Delta_1$ and $\Delta_2$ are defined in Theorem \ref{corAVRIL}.
\item  The map $(c_1,c_2)  \mapsto (c^{+}_1(c_1,c_2), c_2^+(c_1,c_2))$ is continuous.
\end{enumerate}
\end{proposition}

\noindent\emph{Proof of Theorem \ref{EXISTf} assuming Proposition \ref{DANSTHM12}.}\quad
Fix $0<\bar c_1<c_*(f)$  and $0<\epsilon_0<\frac{c_*(f)}{\bar c_1}-1$ small enough so that $Q_{c_1}$ satisfies \eqref{stable}
for all $c_1\in [\bar c_1(1-\epsilon_0),\bar c_1(1+\epsilon_0)]$. Let 
$$\bar c_0 =  \min_{c_1\in [\bar c_1(1-\epsilon_0),\bar c_1(1+\epsilon_0])} c_0(c_1),\quad
\bar K_0= 2 \max_{c_1\in [\bar c_1(1-\epsilon_0),\bar c_1(1+\epsilon_0)]} K_0(c_1),
$$
where $c_0(c_1)$ and $K_0(c_1)$ are defined in Proposition \ref{DANSTHM12}.

Fix an arbitrary $0<\bar c_2<\min(\bar c_0,\epsilon_0^{12})$.
We define $\Omega=[1-\bar c_2^{\frac 1{12}},1+\bar c_2^{\frac 1{12}}]^2$, and the continuous map
$$
\Phi~:~ (\lambda_1,\lambda_2)\in \Omega\mapsto 
\left( \frac {c_1^+(\lambda_1\bar c_1,\lambda_2\bar c_2)}{\bar c_1},\frac {c_2^+(\lambda_1\bar c_1,\lambda_2\bar c_2)}{\bar c_2}\right).
$$
By \eqref{thm12II}, we have
$$
\text{for $j=1,2$,}\quad 
\left|\frac {c_j^+(\lambda_1 \bar c_1,\lambda_2\bar c_2)}{\bar c_j} -\lambda_j \right| \leq \bar K_0 \bar c_2^{\frac {1}{3}}.
$$
This means that 
\begin{equation}\label{TT1}
\|\Phi- \text{Id}\|\leq \bar K_0 \bar c_2^{\frac {1}{3}}.
\end{equation}
Moreover, by possibly taking a smaller $\epsilon_0$,
\begin{equation}\label{TT2}
\text{dist}((1,1),\Phi(\partial\Omega))\geq \bar c_2^{\frac {1}{12}}-\bar K_0 \bar c_2^{\frac {1}{3}}
\geq \frac 12 \bar c_2^{\frac {1}{12}}> \|\Phi-\text{Id}\|.
\end{equation}
From \eqref{TT1} and \eqref{TT2}, we have
$\text{deg}(\Phi,\Omega,(1,1))=\text{deg}(\text{Id},\Omega,(1,1))=1$.
Therefore, from  degree theory there exist $(\bar \lambda_1,\bar \lambda_2)\in \Omega$ such that
$\Phi(\bar \lambda_1,\bar \lambda_2)=(1,1)$ (see for example Theorems 2.3 and 2.1, p30 of  \cite{FG}.)

\smallskip

Now,  for $j=1,2$, we set $c_j=\bar \lambda_j \bar c_j$, and we check that the function 
$u_{c_1,c_2}(t)$ has the property announced in Theorem \ref{EXISTf}. Indeed, since $\Phi(\bar \lambda_1,\bar \lambda_2)=(1,1)$,
we have $c_j^+(c_1,c_2)=\bar c_j$ for $j=1,2$. Moreover, \eqref{thm12I} and \eqref{thm12IV} imply \eqref{TH1A} and \eqref{TH1E}.
Finally, \eqref{TH1B} and \eqref{TH1XX} follow from \eqref{thm12III} and \eqref{defDeltabiscor}.

\medskip

\noindent\emph{Proof of Proposition \ref{DANSTHM12}.}
Let $c_1,c_2$ be as in the statement of Proposition \ref{DANSTHM12} for $0<c_2<c_0(c_1)$ small enough. 
Let  $u(t,x)=u_{c_1,c_2}(t,x)$ be the solution of \eqref{defv12}.
Denote for simplicity $T=T_{c_1,c_2}$ (defined in \eqref{defTccor}).

\smallskip

\noindent\emph{Step 1.}  Control of the modulation parameters of $u(t)$ 
for $t\geq T$.
From Proposition~\ref{INTERACT4col} applied with $T_0=0$ and $\theta=\frac 2{p-1}+\frac 14$,
since $u(0)-v_{c_1,c_2}(0)=0$, we obtain, for some $\rho(t)$,
\begin{equation}\label{atc}
\forall t\in [0,T],\quad \|u(t)-v(t,.-\rho(t))\|_{H^1} \leq K c_2^{\frac 2{p-1}+\frac 14-\frac 1{100}},
\end{equation}
where $|\rho'(t)-c_1|\leq K c_2^{\frac 2{p-1}+\frac 14-\frac 1{100}}$, $ \rho(0)=0$
and so   
\begin{equation}\label{4.7bis}
|\rho(T)-c_1 T|\leq K c_2^{\frac 2{p-1}-\frac 14 -\frac 1{50}}.\end{equation}

By \eqref{3-16cor} and \eqref{atc}, we have
\begin{equation}\label{atcdeux}
 \|u(T)- Q_{c_1}(.-a) - Q_{c_2}(.-b)\|_{H^1}  \leq K c_2^{\frac 2{p-1}+\frac 14-\frac 1{100}},
\end{equation}
for $a=\frac 12 \Delta_1 +\rho(T)$, $b=(c_1-c_2)T + \frac 12 \Delta_2 + \rho(T)$,
so that $$\frac 12 c_1 T\leq a-b\leq 2 c_1 T.$$
Therefore, we can apply Proposition \ref{ASYMPTOTIC} (1) to $u(t)$ with $\omega=\frac 1{p-1}-\frac 1{100}$.
Then, by Claim \ref{LEMMEB1}
we have the decomposition 
of $u(t)$ in terms  of $\eta(t)$, $c_j(t)$, $\rho_j(t)$ ($j=1,2$)
defined for all $t\geq T$:
\begin{equation}\label{4.14B}
\eta(t,x)=u(t,x)- Q_{c_1(t)}(x-\rho_1(t))-Q_{c_2(t)}(x-\rho_2(t)),
\end{equation}
with for  all $t\geq T$, 
\begin{equation}\label{sureta}
\forall t\geq T,\quad
 \|\eta(t)\|_{H^1}\leq K c^{\frac{2}{p-1}-\frac 14 -\frac 1{100}} .
\end{equation}

Now, we claim
\begin{equation}\label{23nov}
|\rho_1(T )-c_1 T -\tfrac 12 \Delta_1 | \leq  K c_2^{\frac 2{p-1}-\frac 12},\quad 
 | \rho_2(T )-c_2 T -\tfrac 12 {\Delta_2}  |\leq K c_2^{\frac 1{p-1}-\frac 1{100}}.
\end{equation}
Proof of \eqref{23nov}.
From \eqref{atc},  \eqref{4.7bis}
and $\|v(T)\|_{H^2}\leq K$, 
we have 
\begin{equation}\label{fdg}
\|u(T )-v(T,.- c_1T )\|_{H^1}   \leq K c_2^{\frac 2{p-1}-\frac 14 -\frac 1{50}}.
\end{equation}
Remark that for $a$ small,
\begin{equation}\label{doublebar}
\tfrac 1K|a|\leq   \|Q_{c_1}-Q_{c_1}(.-a)\|_{L^2}\leq K |a|, 
\quad \tfrac 1K|a|\leq   c_2^{-\frac 1{p-1}+\frac 14} \|Q_{c_2}-Q_{c_2}(.-a)\|_{L^2}\leq K |a|.
\end{equation}
By \eqref{3-16cor} we have
$$
\| v(T ) - Q_{c_1}(.-\tfrac 12 \Delta_1) - Q_{c_2}(.+(c_1-c_2)T  -\tfrac 12 \Delta_2))\|_{H^1}
            \leq K c^{\frac 2{p-1} +\frac 14}.$$
            Thus by \eqref{4.14B}, \eqref{fdg} and \eqref{doublebar}, we deduce \eqref{23nov}.
            
\smallskip

\noindent\emph{Step 2.}  Asymptotic stability.
From   \eqref{sureta}, we can apply  Proposition \ref{ASYMPTOTIC} (2)  
to $u(.+T)$ with $\omega=\frac 1{p-1}-\frac 1{100}$. We deduce that
there exist 
$c_1^+$, $c_2^+>0$,  
 such that
\begin{equation}\label{5-40}
c_j(t)\to c_j^+, \quad 
\rho_j'(t)\to c_j^+,
\quad \text{as $t\to +\infty$,  $j=1,2$,}
\end{equation}
\begin{equation}\label{neuf44}
\lim_{t\rightarrow +\infty}\|w^+(t)\|_{H^1(x> c_2 t/10)}=0,
\end{equation}
where 
$$
w^+(t,x)=u(t,x)-Q_{c_1^+}(x-\rho_1(t))-Q_{c_2^+}(x-\rho_2(t)),
$$
\begin{equation}\label{sept44}
\left|\frac{c_1^+}{c_1}-1\right|\leq   K  c^{\frac 2{p-1}+\frac 14 -\frac 1{100}}  ,\quad
 \left|\frac {c_2^+}{c} - 1\right|\le K c^{\frac 1{p-1}-\frac 1{100}}.
\end{equation}
From \eqref{5-40}, $\|\eta(t)-w^+(t)\|_{H^1}\to 0$
as $t\to +\infty$ and thus, from  \eqref{sureta}, we obtain
$\|w^+(t)\|_{H^1}\leq K  c^{\frac 2{p-1}-\frac 14 -\frac 1{100}}$ for $t$ large.
This concludes the proof of the first part of Proposition \ref{DANSTHM12}.

\medskip

\noindent\emph{Step 3.} Continuity of $c_1^+(c_1,c_2)$ and $c_2^+(c_1,c_2)$.
The proof is the same as in \cite{MMcol1}. Let us give a sketch. 

Let $\bar c_1<c_*(f)$ such that \eqref{stable} holds for $\bar c_1$ and $0<\bar c_2<c_0$ small enough.
First, we prove that the map
$(c_1,c_2)\mapsto c_1^+(c_1,c_2)$ defined in a neighboorhood of $(\bar c_1,\bar c_2)$ 
is continuous.

Denote by
$\eta_{c_1,c_2}(t)$,  $c_{c_1,c_2,j}(t)$, $c_j^+(c_1,c_2)$, the parameters  in the decomposition of   $u_{c_1,c_2}(t)$.
We claim an estimate on $|c_1^+(c_1,c_2)-c_{c_1,c_2,1}(t)| $ which is related to 
the quantities $\mathcal{M}_1(t)$, $\mathcal{E}_1(t)$ defined in section 3.3.
\begin{claim}\label{lllff}
For all $t\geq T_c$,
\begin{equation}\label{unifc}\begin{split}
|c_1^+(c_1,c_2){-}c_{c_1,c_2,1}(t)| &\leq 
  K_0  \int ((\eta_{c_1,c_2})_x^2{+}\eta_{c_1,c_2}^2)(t,x) \psi(x {-} \rho_1(t) {+}c_1 \tfrac {t }4) dx
 +  K_0 e^{-\frac 1 {64} c_1\sqrt{c_2} t }.
\end{split}\end{equation}
\end{claim}

Assuming this claim, let us complete the proof of continuity of $c_1^+(c_1,c_2)$.

Since 
$\|\eta_{\bar c_1,\bar c_2}(t)\|_{H^1(x>\frac {\bar c_2 t}{10})} \to 0$ as $t\to +\infty$, 
for $\varepsilon>0$, there exits
$T_\varepsilon>0$ such that 
 \begin{equation*}\begin{split}
     K_0 \int ((\eta_{\bar c_1,\bar c_2})_x^2+\eta_{\bar c_1,\bar c_2}^2)(T_\varepsilon,x) \psi(x- \rho_1(T_\varepsilon)
     +c_1 \tfrac {T_\varepsilon}4 ) dx
+   K_0  e^{-\frac 1 {64} c_1 \sqrt{c_2 } T_\varepsilon}\leq \varepsilon.
\end{split}\end{equation*}
We fix $T_\varepsilon>0$ to such value. Then, by continuous dependence in $H^1$ of $u_{c_1,c_2}(t)$ solution
of \eqref{kdvf} upon the initial data on $[0,T_\varepsilon]$ (see \cite{KPV})
and of its decomposition in Claim \ref{LEMMEB1}, and the fact that
$u_{c_1,c_2}(0)=v_{c_1,c_2}(0)$ is continuous upon the parameters $(c_1,c_2)$ (see proofs of Proposition \ref{AVRILf} and
Theorem \ref{corAVRIL}), there exists $\delta(\varepsilon)>0$ such that
if $|(c_1,c_2)-(\bar c_1,\bar c_2)|\leq \delta$, then
 \begin{equation*}\begin{split}
     K_0  \int ((\eta_{c_1,c_2})_x^2+\eta_{c_1,c_2}^2)(T_\varepsilon,x) \psi(x- \rho_1(T_\varepsilon)
     +c_1 \tfrac {T_\varepsilon}4 ) dx
+   K_0  e^{-\frac 1 {64}c_1 \sqrt{c_2} T_\varepsilon}\leq 2 \varepsilon,
\end{split}\end{equation*}
\begin{equation*}
|c_{\bar c_1,\bar c_2,1}(T_\varepsilon)-c_{c_1,c_2,1}(T_{\varepsilon})|  \leq  \varepsilon.
\end{equation*}
From Claim \ref{lllff},  applied to $\eta_{c_1,c_2}$, $\eta_{\bar c_1,\bar c_2}$, we have 
$
|c_{1}^+(c_1,c_2)-c_{c_1,c_2,1}(T_\varepsilon)| \leq 2\varepsilon
$ and $|c_{ 1}^+(\bar c_1,\bar c_2)-c_{\bar c_1,\bar c_2,1}(T_\varepsilon)|\leq \varepsilon$. 
Therefore, $
|c_{1}^+(\bar c_1,\bar c_2)-c_{1}^+(c_1,c_2)| \leq 4\varepsilon.
$
Thus, $(c_1,c_2)\mapsto c_{1}^+(c_1,c_2)$ is continuous.

We argue similarly  for $(c_1,c_2)\mapsto c_{2}^+(c_1,c_2)$.
This concludes the proofs of Proposition \ref{DANSTHM12} and of Theorem \ref{EXISTf}.

\medskip

\noindent\emph{Proof of Claim \ref{lllff}.} 
For $T\leq t_0\leq t$, let $\mathcal{M}_1(t)$ and $\mathcal{E}_1(t)$ be defined 
in section 3.3, with $x_0=c_1 \frac {t_0}4$.
From    Claim \ref{NEW} integrated on $[t_0,t]$, 
we obtain
\begin{equation*}\begin{split}
&  \int Q^2_{c_1(t)}-\int Q^2_{c_1(t_0)}   \leq   
( \mathcal{M}_1(t_0)-\mathcal{M}_1(t))
 + K e^{-\frac 1{64} c_1\sqrt{c_2} t_0},
\\ & \left(-E(Q_{c_1(t)})+E(Q_{c_1(t_0)}) -\frac {c_1^+}{100} \Big(\int Q^2_{c_1(t)}-\int Q^2_{c_1(t_0)}\Big)\right)  
\\ & \quad \geq 2\mathcal{E}_1(t)- 2 \mathcal{E}_1(t_0) +\frac 1{100} ( \mathcal{M}_1(t)-\mathcal{M}_1(t_0)) - K  e^{-\frac 1{64} c_1\sqrt{c_2} t_0}.
\end{split}\end{equation*}
Note in particular that 
$
\int_{t_0}^t e^{-\frac 1{16} \sqrt{c_1}(c_1(t-t_0)+x_0)} g_1(t) dt
\leq K e^{-\frac 1 {16} \sqrt{c_1} x_0}\leq K e^{-\frac 1 {64} c_1^{\frac 32} t_0}.
$
Letting $t\to +\infty$, by the asymptotic stability, this gives
\begin{equation*}\begin{split}
 &  \int Q_{c_1^+}^2  -\int Q^2_{c_1(t_0)}  \leq   
 \mathcal{M}_1(t_0) 
 + K e^{-\frac 1{64} \sqrt{c_2} t_0},
\\ &  E(Q_{c_1^+}){-}E(Q_{c_1(t_0)}) {+}\frac {c_1^+}{100} \Big(\int Q^2_{c_1^+} {-}\int Q^2_{c_1(t_0)}\Big) 
 \leq 2 \mathcal{E}_1(t_0) {+} \frac {c_1^+}{100}   \mathcal{M}_1(t_0) 
{+} K e^{-\frac 1{64} c_1 \sqrt{c_2} t_0}.
\end{split}\end{equation*}
By \eqref{quatreet}, we obtain:
\begin{equation*}\begin{split}
|c_1^+-c_1(t_0)| \leq 
K \int (\eta_x^2+\eta^2)(t_0,x) \psi(x - \rho_1(t_0) +\tfrac {t_0}4 ) dx
+K e^{-\frac 1 {64} c_1 \sqrt{c_2} t_0},
\end{split}\end{equation*}
which concludes the proof of Claim \ref{lllff}.

\subsection{Proof of stability. Theorem \ref{STABf}}
 Theorem \ref{STABf} follows directly from Proposition  \ref{INTERACT4col}, Proposition \ref{ASYMPTOTIC}
 and the proof of Theorem~\ref{EXISTf}.
Let $0<\bar c_1<c_*(f)$ such that \eqref{stable} holds for $\bar c_1$. Let $0<\bar c_2<c_0(\bar c_1)$ small enough.
We assume
\begin{equation}\label{hypstab}
\|u(0)-\varphi(0)\|_{H^1}\leq K \bar c_2^{\frac 1{p-1}+\frac 12},
\end{equation}
 where $\varphi=\varphi_{\bar c_1,\bar c_2}$ is the solution constructed  in Theorem \ref{EXISTf}.

From the proof of Theorem \ref{EXISTf}, there exist $(c_1,c_2)$ close to $(\bar c_1,\bar c_2)$ in the following sense (see \eqref{thm12II}):
\begin{equation}\label{proche}
\left|\frac{\bar c_1}{c_1}-1\right|\leq K c_2^{\frac 2{p-1} +\frac 14 -\frac 1{100}},\quad
\left|\frac{\bar c_2}{c_2}-1\right|\leq K c_2^{\frac 1{p-1}-\frac 1{100}},
\end{equation}
so that $\varphi(0)=v_{c_1,c_2}$.
The assumption \eqref{hypstab} on $u(0)$ is thus equivalent to 
\begin{equation}\label{mlml}
\|u(0)-v_{c_1,c_2}(0)\|_{H^1} \leq K c_2^{\frac 1{p-1}+\frac 12}.
\end{equation}
By invariance of \eqref{kdvf} by the transformation $x\to -x$, $t\to -t$, 
it is enough to prove the result for $t\geq 0$.

\medskip

\emph{(i) Estimates on $[0,T_{c_1,c_2}]$.} 

By \eqref{mlml} and
Proposition \ref{INTERACT4col} (applied with $T_0=0$ and $\theta=\frac 1{p-1} +\frac 12$)
we obtain, for all $t\in [0,T_{c_1,c_2}]$, for some $\rho(t)$,
\begin{equation*}
\|u(t)-v(t,x-\rho(t))\|_{H^1} \leq K  c_2^{ \frac 1{p-1}+\frac 12}+K c^{\frac 2{p-1} +\frac 14 -\frac 1{100}}
 \leq K  c_2^{ \frac 1{p-1}+\frac 12},
\end{equation*}
for $c_2$ small.
From \eqref{OUBLI2}, we obtain \eqref{stabbb} on $[0,T_{c_1,c_2}]$.

From Theorem \ref{corAVRIL}, we deduce, for some $a$, $b$, with
$a-b\geq \frac 12 T_{c_1,c_2}$,
\begin{equation}\label{oioi}
\|u(T_{c_1,c_2})-Q_{c_1}(.-a)-Q_{c_2}(.-b)\|_{H^1}\leq K c_2^{\frac 1{p-1}+\frac 12}.
\end{equation}

\medskip

\emph{(ii) Estimates on $[T_{c_1,c_2},+\infty)$.}

 By \eqref{oioi} and Proposition \ref{ASYMPTOTIC}
(applied with $\omega=\frac 14$)
for all $t\in [T_{c_1,c_2},+\infty)$, there exist
$\rho_1(t)$, $\rho_2(t)$
and $c_1^+$, $c_2^+$,  such that 
\begin{equation}\label{triple}\begin{split}
& \|u(t)-Q_{c_1^+}(.-\rho_1(t))-Q_{c_2^+}(.-\rho_2(t))\|_{H^1} \leq K    c_2^{\frac 1{p-1}},\\
& \left|\frac {c_1^+}{c_1}-1 \right| \leq K c_2^{\frac 1{p-1}+\frac 12} , 
\quad  \left|\frac {c_2^+}{c_2} -1 \right| \leq K c^{\frac 14}  .
\end{split}\end{equation}

\medskip

\emph{(iii)} Combining \eqref{proche} and \eqref{triple}, we obtain
$$
\left|\frac {c_1^+}{\bar c_1}-1 \right| \leq K c_2^{\frac 1{p-1}+\frac 12} , 
\quad  \left|\frac {c_2^+}{\bar c_2} -1 \right| \leq K c^{\frac 14}  .
$$
\appendix

\section{Proof of Lemma \ref{surQc2}.}
Proof of \eqref{decay}: it follows from the equation of $Q_c$, \eqref{surf}  and standard arguments.

Note that for any $0<c<c_*$; from \eqref{ellipticf} multiplying by $Q_c'$ and integrating, we get
\begin{equation}\label{surqc}
(Q_c')^2 + 2 F(Q_c) = c\, Q_c^2.
\end{equation}
Using the Taylor decomposition of $F(Q_c)$ (see \eqref{taylor}), we obtain
$$
(Q_c')^2= cQ_c^2 +  \sum_{p+1\leq k_1 \leq k_0 } \sigma_{k_1} Q_c^{k_1}
+O(Q_c^{k_0+1}),
$$
and \eqref{taylor0} follows from $(Q_c^k)'(Q_c^{\widetilde k})'=
k\widetilde k (Q_c')^2 Q_c^{k+\widetilde k-2}.$

Proof of \eqref{taylor1}--\eqref{taylor3}. We prove \eqref{taylor1} and \eqref{taylor3},
\eqref{taylor2} is obtained in a similar way.
Note that from \eqref{ellipticf} and \eqref{taylor}, we get \eqref{taylor1} for $k=1$.
For $k\geq 1$, we have from direct calculations:
\begin{align}
    (Q_c^k)'' 
    & = k(k-1) (Q_c')^2 Q_c^{k-2} + k Q_c'' Q_c^{k-1}\nonumber\\
    & = k(k-1) cQ_c^k - 2k(k-1) Q_c^{k-2} F(Q_c) + ck Q_c^k - k  f(Q_c) Q_c^{k-1} \nonumber \\
    & = k^2 c Q_c^k - 2k(k-1) Q_c^{k-2} F(Q_c) - k f(Q_c) Q_c^{k-1},\label{zz1}
\end{align}
and we get \eqref{taylor1} by using \eqref{taylor} for $f$ and $F$.
Now, we prove  \eqref{taylor3}, from \eqref{zz1},
\begin{align*}
\label{}
    (Q_c^k)^{(4)} 
    & =((Q_c^k)'')''= ck^2 (Q_c^k)'' - 2k(k-1) (Q_c^{k-2} F(Q_c))'' - k (f(Q_c) Q_c^{k-1})''.
\end{align*}
For the first term, we use \eqref{taylor1}. Now, we consider the term $(f(Q_c) Q_c^{k-1})''$,
the term $(Q_c^{k-2} F(Q_c))''$ is similar. We have
\begin{align*}
(f(Q_c) Q_c^{k-1})'' & =
(Q_c^{k-1})'' f(Q_c) + (Q_c')^2 Q_c^{k-2} (2(k-1)   f'(Q_c) + Q_c f''(Q_c)) + Q_c'' Q_c^{k-1} f'(Q_c)
\\ &=
c \left[ (k-1)^2 Q_c^{k-1} f(Q_c) + Q_c^k (2(k-1) f'(Q_c) + Q_c f''(Q_c))+ Q_c^k f'(Q_c)\right]
\\ &
-2F(Q_c)Q_c^{k-2} (2(k-1) f'(Q_c) + Q_c f''(Q_c)) -f(Q_c) Q_c^{k-1} f'(Q_c).
\end{align*}
Now, using Taylor expansions for $f$ (i.e. \eqref{taylor}) and  for
$f'$ and $f''$, we get \eqref{taylor3}. Thus Lemma \ref{surQc2} is proved.

\begin{claim}\label{surDD}
       (i)
        For any integer $r>0$,
    \begin{equation*}
        Q_c^r(y_c) \beta(y_c) = \sum_{\substack{1+r\leq k\leq k_{0}+r \\ 0\leq \ell \leq    \ell_{0}}} 
        c^\ell Q_c^k(y_c) a_{k-r,\ell}.
    \end{equation*}
        (ii) Decomposition of $\beta''$, $\beta^2$, $\beta'\beta$ and $\beta^3$:
    \begin{align*}
        & \beta''(y_c) = \sum_{\substack{1\leq k\leq    k_{0}+p-1 \\ 0\leq \ell \leq    \ell_{0}+1 }}
          c^\ell Q_c^{k}(y_c) {a^{1*}_{k,\ell}}+O(Q_c^{k_0+1}),\quad
        \beta^2(y_c) = \sum_{\substack{2\leq k\leq  2 k_{0} \\ 0\leq \ell \leq 2    \ell_{0} }}
          c^{\ell} Q_c^{k}(y_c)    {a^{2*}_{k,\ell}} ,\\
        &\beta'(y_c)\beta(y_c) = \sum_{\substack{2\leq k\leq    2 k_{0} \\ 0\leq \ell \leq 2    \ell_{0} }} 
          c^{\ell} (Q_c^{k})'(y_c) {a^{3*}_{k,\ell}} ,\quad
        \beta^3(y_c)= \sum_{\substack{3\leq k\leq    3 k_{0} \\ 0\leq \ell \leq 3    \ell_{0} }}
         c^{\ell} Q_c^{k}(y_c)    {a^{4*}_{k,\ell}} ,
    \end{align*}
    where for any $k\geq 1,$ $\ell\geq 0$, the coefficients ${a^{1*}_{k,\ell}}$, ${a^{2*}_{k,\ell}}$, ${a^{3*}_{k,\ell}}$ and ${a^{4*}_{k,\ell}}$ depend on some $(a_{k',\ell'})$ for $(k',\ell')\prec (k,\ell)$.
  \end{claim}
See proof of Claim A.1 in \cite{MMcol1}.
 


\end{document}